\documentclass[a4paper,10pt]{extarticle}
\usepackage{xcolor}
\usepackage[a4paper, total={7in, 8in}]{geometry}
\usepackage{tikz}

\usepackage{cancel}
\usepackage{pgfplots}
\usepackage{amsmath}
\usepackage{mathrsfs}
\usetikzlibrary{arrows}
\usepackage{amsmath}
\usetikzlibrary{calc}
\usepackage{amsmath, amsthm, amsfonts,amssymb,stmaryrd}
\usepackage{hyperref}
\usepackage{pst-node} 
\usepackage[extra]{tipa}
\usepackage{pstricks}
\usepackage{tikz}
\usepackage[spanish,USenglish]{babel}      
\usepackage [T1]{fontenc}
\usepackage{amssymb,amsthm,amsmath,amsfonts}
\usepackage{graphicx}
\usepackage[affil-it]{authblk}
\usepackage{url}
\usepackage{subfigure}
\usepackage{comment}

\usepackage{colortbl} 

\usepackage{graphicx}

\hypersetup{colorlinks=true, linkcolor=blue}

\usepackage{enumerate}
\usetikzlibrary{shapes.geometric}

\newtheorem{thm}{Theorem}[section]
\newtheorem{mthm}{Main Theorem}
\newtheorem{cor}[thm]{Corollary}
\newtheorem{lem}[thm]{Lemma}
\newtheorem{prop}[thm]{Proposition}
\theoremstyle{definition}
\newtheorem{notation}[thm]{Notation}
\newtheorem{defn}[thm]{Definition}
\newtheorem{exam}[thm]{Example}
\theoremstyle{remark}
\newtheorem{rem}[thm]{Remark}

\def\rr{\mathbb{R}}

\newcommand{\coef}[2]{\mathsf{e}_{#1}\left( #2\right) }

\newcommand{\mfp}{\mathfrak{p}}
\newcommand{\mfq}{\mathfrak{q}}

\newcommand{\gdos}{H}

\title{The second order $R$ transform for products of free operators}

\author[1]{Octavio Arizmendi}
\author[2]{Luca Lionni}
\author[3]{Saylé Sigarreta}
\affil[1]{\small Centro de Investigación en Matemáticas, calle Jalisco S/N, Guanajuato, Mexico}
\affil[2]{\small CNRS, ENS de Lyon, LPENSL, UMR5672, 69342 Lyon cedex 07, France}
\affil[3]{ Benemérita Universidad Autónoma de Puebla, Puebla, México  }
\begin{document}
\maketitle

\begin{abstract}
We derive a functional expression for the second-order $R$-transform of the product of two second-order free random variables, expressed in terms of the first- and second-order 
$R$-transforms of the individual factors. The proof is combinatorial and relies on an analysis of the second-order cumulants of the product, together with a recent formula of Muñoz and Perales (2026).

Building on this, we also show in the Appendix how our framework  may be applied to finite free probability. More precisely, we derive an explicit generating series for the infinitesimal part of the finite free cumulants of the finite free multiplicative convolution of two polynomials.
\end{abstract}

\


\section{Introduction}

Free probability, introduced by Voiculescu \cite{VDN92}, has become one of the fundamental tools for studying the asymptotic spectral distributions of large random matrices. A cornerstone of the theory asserts that, under broad assumptions, independent random matrices converge in distribution to freely independent non-commutative random variables as the dimension tends to infinity \cite{V91}. This correspondence provides effective techniques for describing the limiting spectral distributions of various random matrix ensembles and has led to a fruitful interaction between operator algebras and random matrix theory.


To capture not only limiting distributions but also their fluctuations, Mingo and Speicher introduced the notion of second-order freeness \cite{MN04}, extending free probability to a second-order framework. Since then, it has been established that several classical ensembles, including Gaussian and Wishart matrices, are asymptotically second-order free under suitable assumptions \cite{MS06}. As a consequence, second-order free probability has become the natural framework for describing the asymptotic covariances of traces of powers of random matrices, while uncovering deep combinatorial structures through annular non-crossing partitions.

A remarkable feature of free probability is that freeness and second-order freeness admit elegant characterizations in terms of cumulants: free (respectively, second-order free) random variables are precisely those whose mixed first-order (respectively, first and second-order) cumulants vanish \cite{NS06,CMSS07}. This characterization reduces questions about freeness to the analysis of cumulants.

Free cumulants, introduced by Speicher \cite{S94}, arise as the coefficients of Voiculescu's $R$-transform, which plays in free probability a role analogous to that of the logarithm of the characteristic function in classical probability. If $a$ is a self-adjoint non-commutative random variable, its Cauchy transform is defined by
$$
G_a(z)=\varphi({(z-a)^{-1}}),
$$
for $z\in\mathbb C^+$, and the $R$-transform is given by
$$
R_a(z)=G_a^{\langle -1\rangle}(z)-z^{-1},
$$
where $G_a^{\langle -1\rangle}$ denotes the compositional inverse of $G_a$. 

The second-order extension of this theory was developed in \cite{CMSS07}, where second-order free cumulants, the second-order $R$-transform and the combinatorial framework of second-order partitioned permutations were introduced. In particular, the authors established functional relations between first- and second-order moment and cumulant generating functions, providing a fundamental analytic tools for studying fluctuations in free probability.

$$G_u(z,w)= G_u'(z)\,G_u'(w)\,R_u\big(G_u(z),G_u(w)\big)
 +\frac{\partial^{2}}{\partial z\,\partial w}
 \left[
   \log\left(
     \frac{G_u(z)-G_u(w)}{z-w}
   \right)
 \right].$$

According to the above ideas, determining the limiting distribution and fluctuations of a polynomial in freely independent random variables amounts to expressing its first- and second-order free cumulants in terms of the corresponding cumulants of the individual variables, or equivalently, in terms of their first- and second-order $R$-transforms. Understanding how products and sums affect free cumulants is therefore a crucial step toward developing the general theory. The additive case is particularly simple, as the first and second order $R$-transforms were specifically designed to linearize freeness:
$$
R_{X+Y}(z)=R_{X}(z)+R_{Y}(z)\qquad
\text{ and } \qquad
R_{X+Y}(z,w)=R_{X}(z,w)+R_{Y}(z,w).
$$

The multiplicative setting is considerably more involved. Unlike in the classical commutative setting, where multiplication can be reduced to addition through the logarithm, no such reduction is available in the non-commutative framework. Nica and Speicher \cite{NS96} proved that the free cumulants of the product of two free random variables satisfy

$$
\kappa_n(ab)
=
\sum_{\pi\in S_{NC}(n)}
\kappa_\pi^a\,
\kappa_{Kr(\pi)}^b,
$$
revealing the fundamental role played by the Kreweras complement on the lattice of non-crossing partitions. 

This phenomenon persists at second order. More recently, building on the original work of Collins et al. \cite{CMSS07}, a combinatorial formula expressing the second-order cumulants of the product of two second-order free random variables in terms of Kreweras complements on annular partitions was obtained in \cite{GP25}; see Theorem \ref{pt1}.

Motivated by these developments, the main contribution of this paper is to derive a functional expression for the second-order $R$-transform of the product of second-order free random variables. 
Let us present our main results more precisely.

\begin{mthm} []\label{t1}
    Let $a$ and $b$ be two second order free non-commutative random variables, then
\begin{align*}
R_{ab}(x,y)
&= -\frac{\partial^2}{\partial x\,\partial y}
\left[ \log\left(
1
-
x\,y \bigg(\frac{R_a(x\,C_1(x))-R_a(y\,C_1(y))}
         {x\,C_1(x)-y\,C_1(y)}\bigg)\,\bigg(\frac{R_b(x\,C_2(x))-R_b(y\,C_2(y))}
     {x\,C_2(x)-y\,C_2(y)}\bigg)
\right) \right]\nonumber\\[4pt]
&\quad + \,
\big(x\,C_1(x)\big)^{\prime}\,
\big(y\,C_1(y)\big)^{\prime}\,
R_a\big(x\,C_1(x),\,y\,C_1(y)\big)  +
\big(x\,C_2(x)\big)^{\prime}\,
\big(y\,C_2(y)\big)^{\prime}
R_b\big(x\,C_2(x),\,y\,C_2(y)\big).
\end{align*}
\noindent where
\[
C_1(z)=\sum_{n\ge 0} \kappa_{b(ab)^{n}}\, z^n,
\qquad
C_2(z)=\sum_{n\ge 0} \kappa_{a(ba)^{n}}\, z^n,
\]
with
$$\kappa_{b(ab)^{n}}:= \kappa_{n+1}\big(b,ab,\dots,ab\big),\qquad
\kappa_{a(ba)^{n}}:=\kappa_{n+1}\big(a,ba,\dots,ba\big),$$
and for any random variable $u$,

\[R_u(z)=\sum_{n\ge 1}\kappa_n^{\,u}\,z^{\,n-1}~~~~\text{and}  \qquad R_u(z,w)=\sum_{i,j\geq 1} \kappa_{i,j}^u z^{i-1}w^{j-1}.\]

\end{mthm}

To illustrate the scope of our main result, we derive explicit formulas for the second-order cumulants of several natural products of second-order free variables. First, let \(b\) be a second-order free Poisson element with rate \(t>0\), and suppose that \(a\) is second-order free from \(b\). We compute the second-order \(R\)-transform and the second-order cumulants of the product \(ab\) (Proposition \ref{dap}). When \(t=1\), the resulting formula simplifies considerably. In particular, for \(cac^*\), where \(c\) is a second-order circular operator, we recover the result of Arizmendi and Mingo [ \cite{AM23}, Theorem 8.11], which states that conjugation by a circular element transforms second-order free cumulants into second-order moments.

Second, let \(s\) be a second-order standard semicircular operator that is second-order free from \(a\). We derive explicit expressions for the second-order cumulants of \(sas\) in terms of the quantities \(C_j^a(p)\), which are defined through the moments of \(a\) (Proposition \ref{dap1}).

Finally, in the Appendix, we use our main formulas in connection to the recent theory of  finite free probability. Our main application in that setting is the following. Let $\mfp=(p_d)_{d\geq1}$ and $\mfq=(q_d)_{d\geq1}$ be sequences of polynomials with infinitesimal distributions $(\mu,\mu')$ and $(\nu,\nu')$, respectively. Our results then provide an expression for the infinitesimal distribution of their finite free multiplicative convolution $\mfp\boxtimes\mfq=(p_d\boxtimes_d q_d)_{d\geq1}$. More precisely, Theorem~\ref{ffp3} gives the following formula for the generating series of the infinitesimal part of the finite free cumulants of $\mfp\boxtimes\mfq$:
\begin{align*}
\widehat{R}_{\mfp\boxtimes\mfq}(z)
={}& \big(zC_1^{\mu,\nu}(z)\big)'\,
\widehat{R}_{\mfp}\big(zC_1^{\mu,\nu}(z)\big)
+ \big(zC_2^{\mu,\nu}(z)\big)'\,
\widehat{R}_{\mfq}\big(zC_2^{\mu,\nu}(z)\big) \\
&+ \frac{1}{2}\,\frac{d}{dz}
\log\Big(1 - z^2\, R_\mu'\big(zC_1^{\mu,\nu}(z)\big)\,
R_\nu'\big(zC_2^{\mu,\nu}(z)\big)\Big),
\end{align*}
where the series $C_1^{\mu,\nu}$ and $C_2^{\mu,\nu}$ are determined by the system
\[
C_1^{\mu,\nu}(z)=R_{\nu}\big(z\,C_2^{\mu,\nu}(z)\big),
\qquad
C_2^{\mu,\nu}(z)=R_{\mu}\big(z\,C_1^{\mu,\nu}(z)\big).
\]
We apply this to characterize the infinitesimal distributions of various examples including the product of Hermite and Laguerre polynomials.

The paper is organized as follows: Sections \ref{si1} and \ref{si1.5} introduce the necessary preliminaries on non-crossing partitions, permutations, partitioned permutations, first- and second-order probability spaces. Section \ref{si3} presents preliminary results that will be used in the proof of our main theorem. Section \ref{si4} is devoted to the proof of the main result. In Section \ref{si5}, we apply the main theorem to compute the second-order cumulants in several examples. Finally, the Appendix  presents an application to finite free probability, along with the necessary preliminaries. 

\subsubsection*{Acknowledgment}
O. Arizmendi and S. Sigarreta were supported by SECIHTI, Grant number  
CBF-2025-I-4227. L. Lionni was supported by the  ANR JCJC project ``RTFPQuEnt" (ANR-25-CE40-5465).  He also acknowledges support  from the  ANR PRC project  ``TAGADA" (ANR-25-CE40-5672).  We emphasize that no AI of any kind was used during this project.

\section{Non-crossing partitions, permutations and partitioned permutations}\label{si1}

\subsection{The set of partitions and non-crossing permutations}

Consider $M \subset \mathbb{Z}$. A partition of $M$ is a collection of pairwise disjoint subsets $B_i$ of $M$ whose union is $M$. The sets $B_i$ are called \textit{blocks of the partition}. We will denote the set of partitions on  $M$ by $P(M)$. When $M = [n] := \{1,\dots,n\}$, we recover the set of partitions on $n$ elements, simply denoted by $P(n)$. 
 The largest element of $P(M)$, denoted $1_M$, is the partition consisting of a single block.

\paragraph{Non-crossing permutations.}
Consider $M \subset \mathbb{Z}$ and let $\lvert M \rvert$ be its number of elements. A permutation on $M$ is a one-to-one mapping from $M$ to $M$. The set of permutations on $M$ will be denoted by $S_M$, and if $M=[n]$, by $S_n$.  For a permutation $\pi \in S_M$, we will denote by  $\#(\pi)$ the number of cycles in the cycle decomposition of $\pi$. The quantity $|\pi|=|M|-\#(\pi)$ defines a metric on $S_M$. If $\pi, \sigma\in S_M$, we will denote by $\pi\vee\sigma$ the partition of $M$ induced by the transitivity classes of the group generated by $\pi$ and $\tau$ on $M$ (that is, two elements $i,j$ of $M$ are in the same block of $\pi\vee\sigma$ if there exists a finite word $w$ in the alphabet $\{\pi, \sigma\}$, such that $i=w(j)$).

Given two positive integers $n, m$, we introduce the two following elements of $S_n$ and $S_{n+m}$, respectively:
\[
\gamma_n=(1,2\cdots, n),\qquad \textrm{and}\qquad 
\gamma_{m,n}=(1,2,\cdots, m)(m+1, \cdots, m+n).
\]

\begin{defn}{(Non-crossing permutation)}\label{s1d1}
Consider $M\subset \mathbb{Z}$ and a fixed permutation $\gamma\in S_M$. We say that the permutation $\pi \in S_M$ is non-crossing with respect to $\gamma$ if it satisfies the following two conditions:
\begin{enumerate}
    \item $\pi \vee \gamma = 1_M$, 
    \item $\#(\pi)+\#(\pi^{-1}\gamma)+\#(\gamma)=|M|+2$.
\end{enumerate}
The set of non-crossing permutations with respect to $\gamma$ will be denoted by $S_{NC}(\gamma)$. When $M=[n]$ and $\gamma=\gamma_n$, we let $S_{NC}(n):=S_{NC}(\gamma_n)$ (in that case, condition 1 is trivially verified). When $M=[n+m]$ and  
$\gamma = \gamma_{m,n}$, we let $S_{NC}(m,n):=S_{NC}(\gamma_{m,n})$.
\end{defn}

Permutations belonging to  $S_{NC}(n)$ are simply called \emph{non-crossing permutations}. Permutations belonging to $S_{NC}(m,n)$ are called \emph{annular non-crossing permutations. }

\begin{rem}
    The special case $M=[n]$ and $\gamma=\gamma_n$ recovers the classical family of \emph{non-crossing partitions}, denoted by $NC(n)$. We briefly recall this correspondence. A partition $\pi\in P(n)$ is said to be \emph{crossing} if there exist indices $i<j<k<l$ such that $i$ and $k$ belong to one block of $\pi$, while $j$ and $l$ belong to a different block. If no such quadruple exists, then $\pi$ is called \emph{non-crossing}. Once a (cyclic) ordering of each block is fixed, every partition naturally determines a permutation by interpreting its blocks as cycles. A result of Biane \cite{B97} shows that, when the cyclic order is taken to be the natural increasing order, the associated permutation is non-crossing with respect to $\gamma$ in the sense of Definition \ref{s1d1}: this defines a one-to-one correspondence between $S_{NC}(n)$ and $NC(n)$. 
\end{rem}

\begin{figure}
    \centering
    \includegraphics[width=0.1\linewidth]{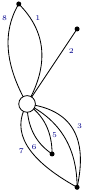}\hspace{2.5cm}    \includegraphics[width=0.3\linewidth]{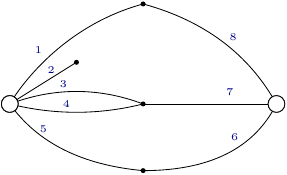}
    \caption{Left: representation of $\pi=(1,8)(2)(3,4,7)(5, 6) \in S_{NC}(8),$ that is, of the combinatorial map $(\gamma_8, \pi)$. Right: representation of $\pi=(1,8)(2)(3,4,7)(5, 6) \in S_{NC}(5,3),$ that is, of the combinatorial map $(\gamma_{5,3}, \pi)$}.
    \label{f1}
\end{figure}

\paragraph{Graphical representations: maps.}
Consider two permutations on a set $M\subset \mathbb{Z}$. From a topological perspective, pairs of permutations $\pi, \sigma\in S_M$ represent labeled bipartite combinatorial maps, which are graphs consisting of vertices partitioned in two sets, respectively colored in black and white, and edges that connect black and white vertices (but are not allowed to connect vertices of the same kind). The edges  are labeled by the elements of $M$, and there is an additional cyclic ordering of the edges attached to every vertex: clockwise for the white vertices, and counterclockwise for the black vertices. Because of the bipartiteness, every element of $M$ is connected a single time to a black (resp.~white) vertex through the edge it labels, so the black and white vertices of a combinatorial map each define a partition of $M$. Moreover, the orderings of edges around each black (resp.~white) vertex defines a cyclic ordering of each block of the aforementioned partition, so that a labeled bipartite combinatorial map thus defines a pair of permutations $(\pi, \sigma)\in S_M \times S_M$.  An edge labeled $i$ is connected to a given cycle (vertex) if $i$ belongs to the corresponding cycle of the associated permutation. If for instance a white (resp.~black) vertex has three edges labeled $i,j,k$ connected to it, and the clockwise (resp.~counterclockwise) ordering of these edges around the vertex is $i$, then $j$, then $k$, then $i$, cyclically, then the first permutation $\pi$ (resp.~the second permutation $\sigma$) has in its cycle decomposition the cycle $(i,j,k)$. Reciprocally, such a pair of permutations defines a labeled bipartite combinatorial map. See the examples in Fig.~\ref{f1}, and \cite{L22} for more details.

The interest of this representation is that the conditions of Def.~\ref{s1d1} have a topological interpretation: the first condition is equivalent to the fact that the graph is connected, while the second condition  is equivalent to the fact that the Euler characteristic of the map is two, that is, that the graph can be drawn on a two-dimensional sphere without crossings, while still respecting the clockwise orderings of edges around the vertices. In such a drawing, the graph delimits some regions of the sphere that are homeomorphic to discs and are in one-to-one correspondence with the cycles in the decomposition of $\pi^{-1}\gamma$ (the same is true considering instead $\gamma\pi^{-1}$, $\gamma^{-1}\pi$, or $\pi\gamma^{-1}$). These regions - called faces - are the 2-cells in Euler's formula, while the 0 and 1-cells are the $\#(\pi) + \#(\sigma)$ vertices and the $n$ edges, respectively. The fact that a map can be drawn on a two-sphere (or the plane) without crossings  explains the term ``non-crossing''. The combinatorial map is also commonly said to be \emph{planar}. For more details, see \cite{W75}.

We are in particular interested in the situation where $\pi\in S_{NC}(\gamma)$, and
\begin{equation}
\label{eq:gamma}
\gamma=(1,\ldots, n_1)\cdots
(n_1+\cdots+n_{m-1}+1, \ldots, n_1+\cdots+n_m),
\end{equation}
for some $m\in \mathbb{N}^\star$ and $n_1, \ldots, n_m\in \mathbb{N}^\star$ satisfying $n_1+\cdots+n_m=n$, where $\mathbb{N}^\star:=\mathbb{N}\setminus \{0\}$. To represent such a $\pi\in S_{NC}(\gamma)$, we will use the map representation of $(\gamma, \pi)$. Hence, the vertices representing the cycles of 
$\gamma$ are thereby represented as white discs.  Since $\gamma$ is given by \eqref{eq:gamma}, the ordering of the labels always grows clockwise around the white discs. The vertices representing the cycles of $\pi$ are represented as black discs. By abuse of language, we will also use the term ``cycles'' to refer to these vertices.  If $\pi\in S_{NC}(m,n)$ (that is, if $\gamma=\gamma_{m,n}$), the cycles of $\pi$ that are connected to edges connected to both cycles of $\gamma_{m,n}$ are said to be \emph{through-cycles}. The edges connected to a through-cycle are called \emph{through-edges}. See Figures \ref{f1}  for concrete examples of these representations.

\begin{figure}
    \centering
    \includegraphics[width=0.45\linewidth]{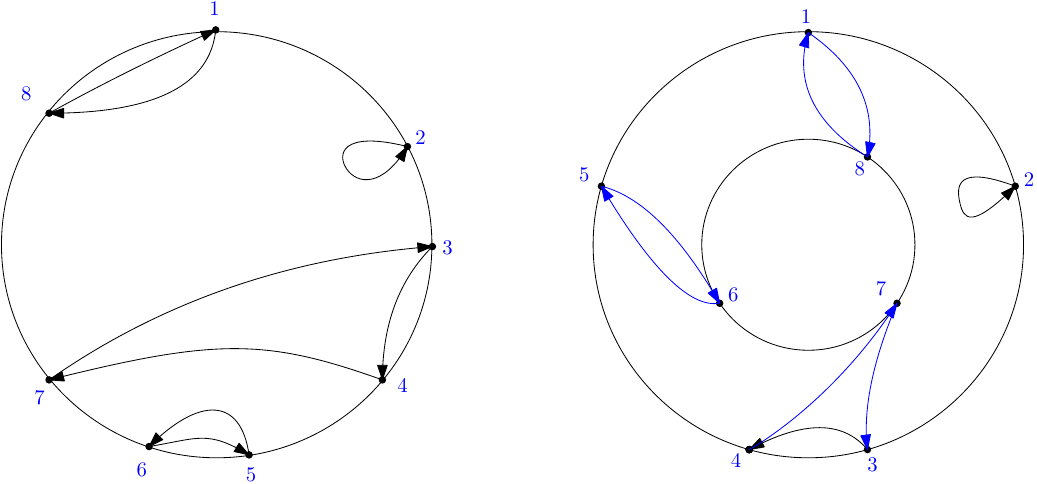}
    \caption{ Left: representation of $\pi=(1,8)(2)(3,4,7)(5, 6) \in S_{NC}(8)$ on the $(8)$-disk. Right: representation of $\pi=(1,8)(2)(3,4,7)(5, 6) \in S_{NC}(5,3)$ on the $(5,3)$-annulus. The through-string are depicted in red.}
    \label{f3}
\end{figure}

There is another common representation of elements of \(S_{NC}(n)\) and \(S_{NC}(m,n)\). The cycles of \(\gamma_n\) and \(\gamma_{m,n}\) are represented by label dots arranged on an \(n\)-disc and an \((m,n)\)-annulus, respectively. The \((m,n)\)-annulus is an annulus in the plane with an outer and an inner boundary circle. The numbers \(1,\ldots,m\) are arranged in clockwise order on the outer circle, while the numbers \(m+1,\ldots,m+n\) are arranged in counterclockwise order on the inner circle. For a permutation \(\pi\), if \(\pi(i)=j\), we draw a \emph{string} connecting the dots labeled \(i\) and \(j\). Under this representation, $\pi\in S_{NC}(n)$ and $\pi\in S_{NC}(m,n)$ corresponds to a collection of cycles that can be drawn without crossings in the corresponding $n$-disc and $(m,n)$-annulus. Each cycle encloses a region homeomorphic to a disk, with the boundary of the region oriented clockwise. In the annular case, at least one cycle connects the inner and outer boundary circles. The two representations are readily seen to be in bijection. For \(\pi\in S_{NC}(m,n)\), the term ``through-string'' is commonly used in the literature for a string connecting the inner and outer boundary circles. Thus, \emph{through-edges} and \emph{through-strings} do not refer to the same concept. See Figure~\ref{f3} for concrete examples of this representation.

\paragraph{Kreweras complement.}
Given $\gamma\in S_n$ and $\pi\in S_{NC}(\gamma)$, the Kreweras
complement of $\pi$, denoted by $Kr(\pi)$, is defined by
\begin{equation}
\label{eq:def-kreweras}
Kr(\pi):=\pi^{-1}\gamma.
\end{equation}

To represent both permutations in a unified way, we regard $Kr(\pi)$
as acting on the barred set $\{\bar{1},\ldots,\bar{n}\}$ via the
natural identification $i\leftrightarrow\bar{i}$. With this convention, the permutation $\pi\times Kr(\pi)$ acts on the set 
$$\operatorname{Int}[n]:=\{1, \bar 1, \ldots, n, \bar n\} \quad \text{as} \quad
\pi\times Kr(\pi) :\ \left\{ \begin{array}{ll}
    i & \mapsto\ \pi(i) \\
     \bar i & \mapsto\ Kr(\pi)(\bar i)
\end{array} \right. .$$

Given a permutation $\gamma\in S_n$, we denote by
$\operatorname{Int}(\gamma)$ the permutation on $\operatorname{Int}[n]$ defined by 
\[
\operatorname{Int}(\gamma):\ \left\{ \begin{array}{ll}
    i & \mapsto\ \bar i \\
     \bar i & \mapsto\ \gamma(i)
\end{array} \right. .
\]
Thus, $\operatorname{Int}(\gamma)$ preserves the cycle structure of
$\gamma$ while interlacing each element with its barred counterpart. If for instance $\gamma$ is given by \eqref{eq:gamma}, then:
\[
\operatorname{Int}(\gamma)
=
(1,\bar{1},\ldots,n_1,\overline{n_1})
\cdots
(n_1+\cdots+n_{m-1}+1,\overline{n_1+\cdots+n_{m-1}+1},
\ldots,n_1+\cdots+n_m,\overline{n_1+\cdots+n_m}).
\]

We also define the permutation $\nu=\nu(\gamma, \pi)$ on $\operatorname{Int}[n]$ as:
\begin{equation}
\label{eq:def-nu}
\nu = \bigl(\pi\times Kr(\pi)\bigr)^{-1}\operatorname{Int}(\gamma). 
\end{equation}

We will need the following properties of  $Kr(\pi)$.

\begin{lem}
\label{lem:Kreweras-prop}
    Consider $\pi, \gamma\in S_n$ such that $\pi\in S_{NC}(\gamma)$. Then the following properties are true:
    \begin{enumerate}
    \item $\nu(\gamma, \pi)$ is an involution without fixed points.
    \item $\pi\times Kr(\pi)\in S_{NC}(\operatorname{Int}(\gamma))$. 
    \end{enumerate}
\end{lem}
Because of the definition \eqref{eq:def-nu}, the second statement implies that $\nu(\gamma, \pi)=Kr(\pi\times Kr(\pi))$.

\begin{proof} Regarding the first statement, by definition, we have
\begin{equation*} \nu(\gamma, \pi):\ \left\{ \begin{array}{ll} i & \mapsto\ Kr(\pi)^{-1}(\bar i) \\ \bar i & \mapsto\ \pi^{-1}\gamma(i)=\overline{Kr(\pi)(\bar i)}\end{array} \right. . \end{equation*}
The last equality holds under the convention that $\bar{\bar{j}}=j$ for $j\in [n]$.
Consequently, $\nu\neq id_{\operatorname{Int}[n]}$ and $\nu^2=\operatorname{id}_{\operatorname{Int}[n]}.
$ In particular,$$
\#\bigl(\nu(\gamma,\pi)\bigr)=n.
$$ Since $\pi\in S_{NC}(\gamma)$, the first condition in Definition~\ref{s1d1} directly implies that

$$
\pi\times Kr(\pi)\vee\operatorname{Int}(\gamma)
=1_{\operatorname{Int}[n]}.
$$ Moreover, by construction,

$$
\#\bigl(\operatorname{Int}(\gamma)\bigr)=\#(\gamma)
\quad \text{and} \quad 
\#\bigl(\pi\times Kr(\pi)\bigr)
=\#(\pi)+\#\bigl(\pi^{-1}\gamma\bigr).
$$ Therefore, since $\pi\in S_{NC}(\gamma)$, by the second condition in Definition~\ref{s1d1}, we obtain

$$
\#\bigl(\pi\times Kr(\pi)\bigr)
+\#\bigl(\nu(\gamma,\pi)\bigr)
+\#\bigl(\operatorname{Int}(\gamma)\bigr)
=2n+2.
$$ This completes the proof.

\end{proof}

In fact,  given $\pi\in S_{NC}(\gamma)$, it can be shown that there is a unique permutation $\sigma\in S_n$ such that $\pi\times \sigma\in S_{NC}(\operatorname{Int}(\gamma))$ and $(\pi\times \sigma)^{-1}\operatorname{Int}(\gamma)$ is an involution without fixed point, so this provides an alternative definition of  the Kreweras complement \eqref{eq:def-kreweras}. It is equivalent to saying that $Kr(\pi)$ is the permutation $\sigma$ on $\{\bar 1, \ldots, \bar n\}$ with the largest possible cycles that satisfies $\pi\times \sigma\in S_{NC}(\operatorname{Int}(\gamma))$.  In the case where $\gamma=\gamma_n$ for instance, and under the above bijection between non-crossing partitions and
non-crossing permutations, this can be reformulated as saying that given $\pi\in NC(n)$, $Kr(\pi)$ is the coarsest non-crossing partition on $\{\bar 1, \ldots, \bar n\}$ such that $\pi\times Kr(\pi)$ is non-crossing on $\{1,\bar{1}, 2, \bar 2, \ldots, n,\bar{n}\}$ (ordered from left to right), a definition commonly used in the literature.  

\begin{figure}[h!]
    \centering
    \includegraphics[width=0.3\linewidth]{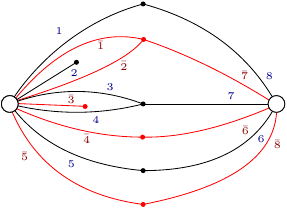}\hspace{2.5cm}\includegraphics[width=0.3\linewidth]{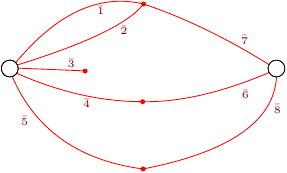}
    \caption{Map representation of $\bigl(\operatorname{Int}(\gamma_{5,3}), \pi\times Kr(\pi)\bigr)$ (left) and $Kr(\pi)$ (right), where $\pi=(1,8)(2)(3,4,7)(5, 6) \in S_{NC}(5,3),$ is represented in Fig.~\ref{f1}.}
    \label{f2}
\end{figure}

Given $\pi\in S_{NC}(\gamma)$, one  can use these properties to easily construct $Kr(\pi)$ in the map representation, as illustrated in  Fig.~\ref{f2}. The recipe is the following. Consider a drawing of the map representing $(\gamma, \pi)$ on the plane, without crossings. 
\begin{itemize}
\item  In each face (the regions on the plane delimited by the graph), draw a red disc.
Add, interleaved with the existing edge labels, the barred counterpart label to the boundary of the white discs. 
\item 
For each face, add red edges connecting the red disc to all the barred counterpart labels that bound the face, in a non-crossing matter.
\end{itemize}
At this point, one has drawn the combinatorial map representing $\bigl(\operatorname{Int}(\gamma), \pi\times Kr(\pi)\bigr)$. In particular, by construction it has no crossings, which is equivalent to the second statement of Lemma~\ref{lem:Kreweras-prop}. It also has the property that every face is bounded by four edges (every face is a square), with repetition (that is, if a red disc has a single edge attached, it is counted twice). This is the first statement of Lemma~\ref{lem:Kreweras-prop}, i.e.~$\nu(\gamma, \pi)$ defined in \eqref{eq:def-nu} is an involution without fixed points\footnote{This is easily seen graphically for instance: one has that $\nu(\gamma, \pi) = \bigl(\pi\times Kr(\pi)\bigr)^{-1} \operatorname{Int}(\gamma) = Kr(\pi\times Kr(\pi)) $. Applying the procedure to construct the Kreweras complement $\nu(\gamma, \pi)$ of $\pi\times Kr(\pi) $, since the faces are squares, every vertex added in this procedure will have two edges connected to it, so the cycle decomposition of $\nu(\gamma, \pi)$ only consists of transpositions. Reciprocally, for the cycles of $\nu(\gamma, \pi)$ to be all transpositions, the faces of  $ \pi\times Kr(\pi)$ must be squares. }. If we now remove the black discs and edges to keep only the white discs and the red discs and edges, we obtain $Kr(\pi)$.

\section{First- and second-order free probability}\label{si1.5}

A \emph{(first-order) non-commutative probability space} is a pair
$(\mathcal A,\varphi)$, where $\mathcal A$ is a unital algebra over
$\mathbb C$ and $\varphi:\mathcal A\to\mathbb C$ is a linear functional satisfying
$\varphi(1)=1$. Let $(\mathcal A,\varphi)$ be a non-commutative probability space, and let
\[
\kappa_n:\mathcal A^{\otimes n}\to\mathbb C,
\qquad n\ge1,
\]
be a family of multilinear functionals. We extend this multiplicatively to a family of multilinear functionals, called the free cumulants
\[
\kappa_\pi:\mathcal A^{\otimes n}\to\mathbb C,\quad n \geq 1,\quad \pi\in S_{NC}(n).
\]  
by
\[
\kappa_\pi(a_1,\ldots,a_n)
=
\prod_{\substack{B\in\pi\\ B=(i_1,\ldots,i_l)}}
\kappa_l(a_{i_1},\ldots,a_{i_l}),
\]
where the indices in each cycle satisfy
\[
i_1<i_2<\cdots<i_l.
\] This construction can then be used to define the sequence
$\{\kappa_n\}_{n\geq 1}$
through the moment-cumulant formula:
\begin{equation}
\label{eq:moment-free-cumulant-1st-order}
\varphi(a_1\cdots a_n)
=
\sum_{\pi\in S_{NC}(n)}
\kappa_\pi^{\vec a},
\end{equation}
where $\vec a = (a_1, \ldots, a_{n})$.

A \emph{second-order non-commutative probability space} is a triple
$(\mathcal A,\varphi,\varphi_2)$, where $(\mathcal A,\varphi)$ is a tracial
non-commutative probability space, that is,
\[
\varphi(ab)=\varphi(ba),
\qquad a,b\in\mathcal A,
\quad \text{and} \quad
\varphi_2:\mathcal A\times\mathcal A\to\mathbb C
\]
is a bilinear, symmetric functional that is tracial in each variable and satisfies
\[
\varphi_2(1,a)=\varphi_2(a,1)=0,
\qquad a\in\mathcal A.
\]

Analogously, the second-order free cumulants,

\[
\kappa_{m,n}
(a_1,\ldots,a_{m};
a_{m+1},\ldots,a_{m+n}
),
\]
 are characterized by the second-order moment--cumulant formula \cite{CMSS07}.

\begin{equation}
\label{comb2}
\varphi_2(a_1\cdots a_m,\,
a_{m+1}\cdots a_{m+n})
=
\sum_{\pi\in S_{NC}(m,n)} \kappa_{\pi}^{\vec a}
+ 
\sum_{\pi=\pi_1\times \pi_2\in S_{NC}(p)\times S_{NC}(q)}
\ \sum_{\substack{U_1\in\pi_1\\[1pt] U_2\in\pi_2}}
    \kappa^{{\vec{a}}}_{|U_1|,\,|U_2|}\,
    \kappa^{{\vec a}}_{\pi\setminus\{U_1,U_2\}}, 
\end{equation}
where $\vec a = (a_1, \ldots, a_{m+n})$.

Let $u$ be a random variable in a first-order non-commutative probability space $(\mathcal A,\varphi)$. Throughout the paper, we use the notation
\[
\alpha_n^u=\varphi(u^n),
\qquad
\kappa_n^u=\kappa_n(u,\ldots,u),
\qquad n\ge 1.
\]

Similarly, if $(\mathcal A,\varphi,\varphi_2)$ is a second-order non-commutative probability space, we write
\[
\alpha_{n,m}^u=\varphi_2(u^n,u^m),
\qquad
\kappa_{n,m}^u=\kappa_{n,m}(u,\ldots,u;\,u,\ldots,u),
\qquad n,m\ge 1.
\]

For a non-commutative random variable \(u\), the combinatorial relation between
moments and cumulants in first order can be expressed in terms of generating
functions. The Cauchy transform is
\[
G_u(z)=\frac{1}{z}+\sum_{n\ge 1}\alpha_n^{\,u}\,z^{-(n+1)},
\]
and the \(R\)-transform is
\begin{equation}
\label{def:R-transform}
    R_u(z)=\sum_{n\ge 1}\kappa_n^{\,u}\,z^{\,n-1}.
\end{equation}

The combinatorial relation \eqref{eq:moment-free-cumulant-1st-order} can be reformulated in terms of the generating functions above as:
\begin{equation}\label{eq1}
    \frac{1}{G_u(z)}+R_u(G_u(z))=z.
\end{equation}  
It may also be formulated in terms of the generating functions
\[
M_u(z)=1+\sum_{n\geq 1} \alpha_n^{u}z^n, \qquad C_u(z)=1+\sum_{n\ge 1}\kappa_n^{\,u}\,z^{\,n},
\]
as:
\begin{equation}\label{eq2}
  G_u(z)=\frac{ M_u(1/z)}{z}, \qquad  R_u(z)=\frac{C_u(z)-1}{z}.
\end{equation}  and

\begin{equation}\label{eq21}
M_u(z)=C_u(z\,M_u(z)).
\end{equation}  For a more detailed discussion, we refer the reader to \cite{NS06}.

The $S$-transform of $u$ is defined by
\[
S_u(z) := \frac{1+z}{z}\, (M_u-1)^{\langle -1\rangle}(z),
\] In particular, replacing $z$ with $(M_u-1)^{\langle -1\rangle}(z)$ into 
\eqref{eq21}, we have that $(C_u(z)-1)^{{\langle -1\rangle}}=zS_u(z)$. Hence, in view of \eqref{eq2}, 

\begin{equation}\label{eq22}
R_u(z)\, S_u\left(z\,R_u(z)\right) = 1.
\end{equation} Note that, a priori, the $S$-transform of $u$ is defined only when 
$\alpha^u_1 \neq 0$. However, by allowing formal power series in $\sqrt{z}$, 
the definition of the $S$-transform is viable when 
$\alpha^u_1 = 0$ and $ \alpha^u_2 \neq 0$. In this setting, the previous formulas
remains valid (see~\cite{RS07}).  The importance of the $S$-transform in free probability theory stems from the
fact that it provides a suitable machinery for studying products of free random
variables. This idea is formally captured in the following result, which
constitutes an extended version of the original statement (see Lecture 18 from \cite{NS06} and Theorem 2.5. from ~\cite{RS07}). Let $a$ and $b$
be free random variables such that $\alpha^b_1 \neq 0$ and either  $\alpha^a_1 \neq 0$ or $\alpha^a_1 = 0$  and $\alpha^a_2 \neq 0$. Then
\[
S_{ab}(z) = S_a(z)\, S_b(z).
\] The second order analogues of $G_u(z)$ and $R_u(z)$ are defined by
\[
G_u(z,w)=\sum_{n,m\ge 1}\alpha_{n,m}^{\,u}\,z^{-(n+1)}\,w^{-(m+1)}
\quad \text{and} \quad
R_u(z,w)=\sum_{n,m\ge 1}\kappa_{n,m}^{\,u}\,z^{\,n-1}w^{\,m-1}.
\]
The combinatorial relations \eqref{comb2} can be reformulated in terms of these generating function as \cite{CMSS07}:
\begin{equation}\label{eq3}
\begin{aligned}
G_u(z,w)
&=\frac{\partial^{2}}{\partial z\,\partial w}
 \left[
   \log\left(
     \frac{G_u(z)-G_u(w)}{z-w}
   \right)
 \right] +  G_u'(z)\,G_u'(w)\,R_u\big(G_u(z),G_u(w)\big).
\end{aligned}
\end{equation} Also, we consider the auxiliary series
\[
\widehat C_u(z,w)
=
\sum_{i,j\ge 0}
\kappa_{i+j+2}^{\,u}\,
z^{\,i+1}w^{\,j+1},
\]
which satisfies
\begin{equation}\label{ecr}
\widehat C_u(z,w)
=
z\,w \frac{R_u(z)-R_u(w)}{z-w}.\end{equation} 
For a more detailed account of the second-order construction introduced above, we refer the reader to \cite{CMSS07}.

\section{Preliminary Results for the Proof of the Main Theorem}\label{si3}

Given $a$ and $b$  two second order free non-commutative random variables, in addition to the moment and cumulant series of $a$, $b$ and $ab$, we consider also  moments series of the form

\[
C_1(z)=\sum_{n\ge 0} \kappa_{b(ab)^{n}}\, z^n,
\qquad
C_2(z)=\sum_{n\ge 0} \kappa_{a(ba)^{n}}\, z^n,
\]
with
$$\kappa_{b(ab)^{n}}:= \kappa_{n+1}\big(b,ab,\dots,ab\big),\qquad
\kappa_{a(ba)^{n}}:=\kappa_{n+1}\big(a,ba,\dots,ba\big).$$

We will use the following elementary properties of free cumulants, which follow from the product formula for cumulants; see Theorems 11.12 and 14.4 in \cite{NS06}.

\begin{prop} \label{p1}
Let $a$ and $b$ be two  non-commutative random variables.
The following properties  hold for $n\geq 0$:

\medskip
\noindent\textbf{a)} 
    \[
    \kappa_{n+1}(ab,\ldots,ab)
    =
    \kappa_{n+2}(a,b,ab,\ldots,ab)
    +
    \sum_{i=1}^{n+1}
        \kappa_{i}(b,ab,\ldots,ab)\,
        \kappa_{\,n+2-i}(a,ba,\ldots,ba).
    \]

\medskip
\noindent Moreover, if $a$ and $b$ are free, then

   \medskip
\noindent\textbf{b)} \[\  \kappa_{n+1}(ab,\ldots,ab)=\kappa_{n+1}(ba,\ldots,ba).\]

 \medskip
\noindent\textbf{c)} \[\  \kappa_{n+1}(ab,\ldots,ab,a)
    =\kappa_{n+1}(a,ba,\ldots,ba).\]

\noindent\textbf{d)}  
    \[
        \kappa_{n+2}(a,b,ab,\ldots,ab)=0.
    \]

\medskip
\noindent\textbf{e)}  $$
    \kappa_{n+1}(b,ab,\ldots,ab)= \sum_{\pi \in S_{NC}(n)}\kappa_{\pi}^a\kappa^b_{Kr^{+1}(\pi)},$$ where $Kr^{+1}(\pi)$
denotes the permutation obtained from  $Kr(\pi)$ by inserting the element $\bar 0$ at the beginning of the cycle whose last element is $\bar n$.

\medskip
\noindent\textbf{f)} $$
    \kappa_{n+1}(ab,\ldots,ab,a)= \sum_{\substack{\pi \in {S_{NC}(n+1)} \\
    \pi(n+1)=1
   }
    } \kappa_{\pi}^a\kappa^b_{Kr^{-1}(\pi)},$$ where  $Kr^{-1}(\pi)=Kr(\pi) \setminus(\overline{n+1}).$
\end{prop}

\begin{prop} \label{p2} 
Let $a$ and $b$ be two free non-commutative random variables.  Given

\[
C_1(z)=\sum_{n\ge 0} \kappa_{b(ab)^{n}}\, z^n,
\qquad
C_2(z)=\sum_{n\ge 0} \kappa_{a(ba)^{n}}\, z^n,
\]
with
$$\kappa_{b(ab)^{n}}= \kappa_{n+1}\big(b,ab,\dots,ab\big),\qquad
\kappa_{a(ba)^{n}}=\kappa_{n+1}\big(a,ba,\dots,ba\big),$$
and for any random variable $u$, 

\[C_u(z)=1+\sum_{n\ge 1}\kappa_n^{u}\,z^n,  \qquad M_u(z)=1+\sum_{n\geq 1} \alpha_n^{u}z^n,\] it follows that

\medskip
\noindent\textbf{a)} $R_{ab}(z)=R_{ba}(z)=C_1(z) \cdot C_2(z).$

\medskip
\noindent\textbf{b)}  $M_{ab}(z)=M_{ba}(z)=C_b(z\cdot M_{ab}(z)\cdot C_2(z\,M_{ab}(z))).$

\medskip
\noindent\textbf{c)} $M_{ba}(z)=M_{ab}(z)=C_a(z\cdot M_{ab}(z)\cdot C_1(z\,M_{ab}(z))).$

\medskip
Moreover, if $\alpha_1^b \neq 0$ and either $\alpha_1^a \neq 0$ or $\alpha_1^a = 0$ and $\alpha_2^a \neq 0$, then

\medskip
\noindent\textbf{d)} $C_{ab}(z)=C_{ba}(z)=C_{a}(z\,C_1(z))=C_{b}(z\,C_2(z)).$

\medskip
\noindent\textbf{e)} $C_2(z)=R_{a}(z\,C_1(z))$ and $C_1(z)=R_{b}(z\,C_2(z)).$

\end{prop}

\begin{proof}
The first statement follows directly from Proposition~\ref{p1}\, $a)$, $b)$ and $d)$.

For the second statement, by the moment-cumulant formula \eqref{eq:moment-free-cumulant-1st-order}, we have
\[
\varphi\bigl(a(ba)^n\bigr)
=
\sum_{\pi \in {S_{NC}(n+1)}}
\kappa_{\pi}\bigl(a,ba,\ldots,ba\bigr).
\] Given $\pi \in S_{NC}(n+1)$, let  $V$ denote the cycle of $\pi$ that contains the element $1$ (which corresponds to the only argument which is $a$ and not $ab$).
Then $\pi$ admits a decomposition of the form
\[
\pi = V \times \pi_1 \times \cdots \times \pi_r,
\quad \text{where} \quad
V = (1,i_1 + 2, i_1 + i_2 + 3, \ldots, i_1 + \cdots + i_{r-1} + r),
\]
for some  $1 \le r \le n+1$, 
$i_\ell \in \{0,1,\ldots,n+1-r\}$  such that $\sum_{\ell=1}^ri_\ell+r=n+1$ and 
$\pi_\ell \in S_{NC}(i_\ell)$ for $\ell = 1,\ldots,r$.
Consequently, applying this decomposition to all $\pi \in S_{NC}(n+1)$ , and performing independently the and sums over the $\pi_\ell$ using \eqref{eq:moment-free-cumulant-1st-order} yields
\[
\varphi\bigl(a(ba)^n\bigr)
=\sum_{r=1}^{n+1}\sum_{\substack{ \{i_{\ell}\}_1^r \\
\sum_{\ell=1}^r i_\ell+r=n+1}}
\kappa_r(a,ba,\dots,ba) \prod_{\ell=1}^r\varphi\bigl((ba)^{i_\ell}\bigr) . 
\] Hence,

\[
M_2(z)
:=
\sum_{n \geq 0} \varphi\bigl(a(ba)^n\bigr) z^n
=
M_{ba}(z)\, C_2\bigl(z M_{ba}(z)\bigr).
\] On the other hand, again by the moment-cumulant formula \eqref{eq:moment-free-cumulant-1st-order},
\[
\varphi\bigl((ba)^n\bigr)
=
\sum_{\pi \in S_{NC}(2n)}
\kappa_{\pi}\bigl(b,a,\ldots,b,a\bigr).
\]
Similarly, using the block containing \(1\) as a reference to describe \(\pi\) and from the vanishing of mixed cumulants, we obtain
\begin{equation} \label{eq0}
    M_{ba}(z)
=
C_b\bigl(z M_2(z)\bigr).
\end{equation}
Combining the two identities and Proposition~\ref{p1}\ $b)$ yields statement $b)$.
The proof of~$c)$ follows in the same manner.

Let us address part $d)$. Substituting $z$ by $(M_{ab}-1)^{\langle -1\rangle}(z)$ 
into $c)$ we obtain
\[
1+z=C_a\bigl(z\,S_{ab}(z)\,C_1(z\,S_{ab}(z))\bigr).
\]
Hence, replacing $z$ with $z\,R_{ab}(z)$, according to
\eqref{eq22} and \eqref{eq2}, yields the desired result. The remaining part
follows straightforwardly. 

Now consider $e)$. According to \eqref{eq2} and $d)$, it follows that
\[
R_{ab}(z)=C_1(z)\,R_a(z\,C_1(z)).
\]
Finally, by part $a)$, the first part is complete. The second part of the statement follows in a  similar way.

\end{proof}

\begin{rem}The generating series \(M_1(z)\) and \(M_2(z)\) were originally defined in \cite{RS07}. Moreover, the  functional relation of Equation~(\ref{eq0})  corresponds to Equation~(1) in \cite{RS07}.

\end{rem}

\section{Proof of the Main Theorem} \label{si4}
In this section, we prove the main theorem. Our starting point is a reformulation of Theorem 7.9 from \cite{CMSS07} for the case when $a_n=a$ and $b_n=b$ for all $n \geq 1$, as recently presented by Munoz and Perales \cite{GP25}. This alternative version is particularly useful for our purposes, as it provides a sum indexed by non-crossing permutations and their Kreweras complements.

\begin{thm} (Theorem 1.7. \cite{GP25}) \label{pt1} Let $a,b$ be two second-order free random variables. Then for any $p,q\ge 1$,

\begin{equation*}
\begin{aligned}
\kappa_{p,q}^{ab}
={}&\sum_{\pi\in S_{NC}(p,q)}
\kappa^{a}_{\pi}\kappa^{b}_{Kr(\pi)}
+\sum_{\pi=\pi_1\times\pi_2\in S_{NC}(p)\times S_{NC}(q)}
\sum_{U_1\in\pi_1,\;U_2\in\pi_2}
\kappa^{a}_{|U_1|,|U_2|}
\kappa^{a}_{\pi\setminus\{U_1,U_2\}}
\kappa^{b}_{Kr(\pi_1)\times Kr(\pi_2)}
\\[-1pt]
&+\sum_{\pi=\pi_1\times\pi_2\in S_{NC}(p)\times S_{NC}(q)}
\sum_{U_1\in Kr(\pi_1),\;U_2\in Kr(\pi_2)}
\kappa^{b}_{|U_1|,|U_2|}
\kappa^{b}_{Kr(\pi_1)\times Kr(\pi_2)\setminus\{U_1,U_2\}}
\kappa^{a}_{\pi}.
\end{aligned}
\end{equation*}

\end{thm}

The proof of the Main Theorem will be split in two lemmas, considering separately each summand in the formula from the above theorem. They follow the same logic as Propositions 6.1 and 6.2 of \cite{CMSS07}, but the decompositions need to be adapted to include the contributions of the Kreweras. 

We remind the definitions
\begin{equation}
\label{eq:C1C2}
    C_1(z)=\sum_{n\ge 0} \kappa_{b(ab)^{n}}\, z^n,
\qquad
C_2(z)=\sum_{n\ge 0} \kappa_{a(ba)^{n}}\, z^n,
\end{equation}

We will refer to the cycle of $\gamma_{p,q}$ containing 1 as its ``first cycle''.

\begin{lem} []\label{l1}
    Let $a$ and $b$ be two second-order free random variables, $p,q$ two positive integers, and let
\[
\kappa^{ab,1}_{p,q}=\sum_{\pi\in S_{NC}(p,q)}
    \kappa^{a}_{\pi}
    \kappa^{b}_{K r(\pi)}.
\]
 Then the series
\[
R_{ab,1}(x,y)
:=\sum_{p,q\ge 1}\kappa^{ab,1}_{p,q}\, x^{p-1}y^{q-1},
\]
can be expressed as

\begin{align*}
R_{ab,1}(x,y)
&= - \frac{\partial^2}{\partial y \,\partial x}
\left[ \log\left(
1
-
x\,y \bigg(\frac{R_a(x\,C_1(x))-R_a(y\,C_1(y))}
         {x\,C_1(x)-y\,C_1(y)}\bigg)\,\bigg(\frac{R_b(x\,C_2(x))-R_b(y\,C_2(y))}
     {x\,C_2(x)-y\,C_2(y)}\bigg)
\right) \right]\nonumber.\\
\end{align*}

\end{lem}

\begin{figure}[ht]
    \centering
    \includegraphics[width=0.65\textwidth]{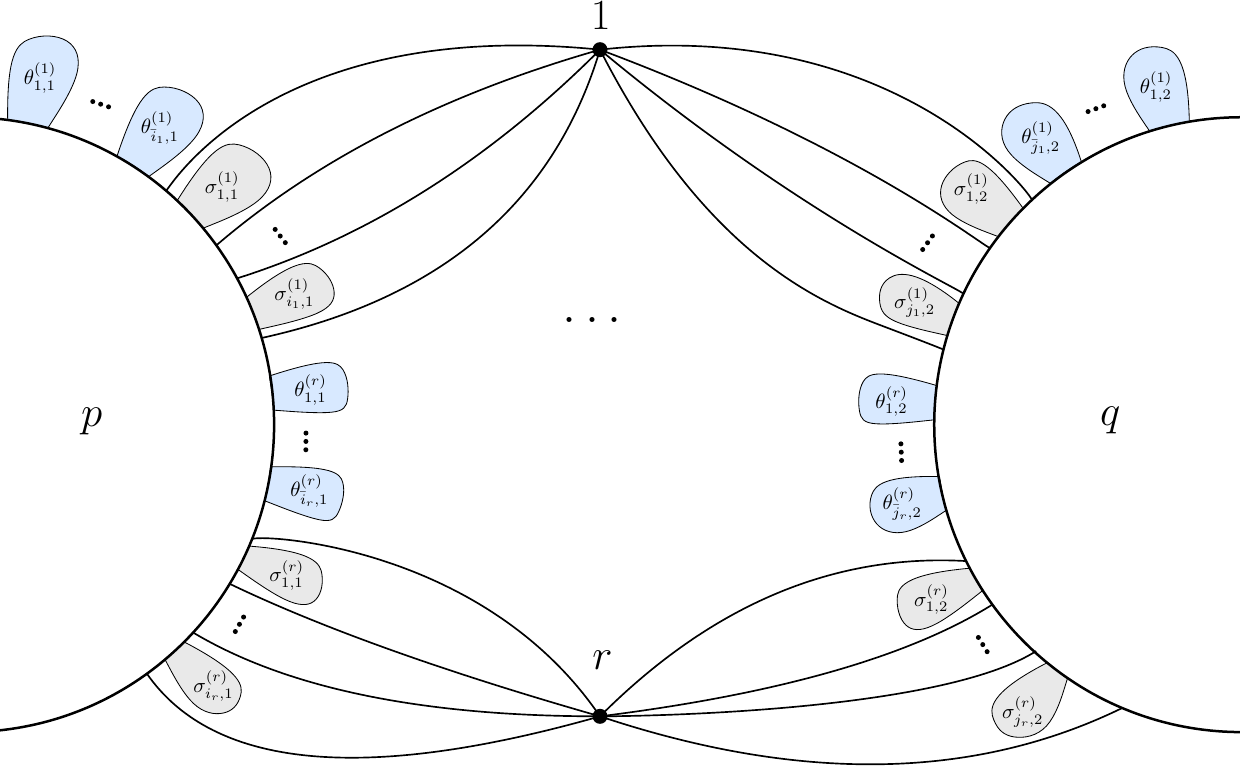}
    \caption{Representation of $\pi \in S_{NC}(p,q)$ appearing in Lemma \ref{l1}.}
    \label{f4}
\end{figure}

\begin{figure}[ht]
    \centering
    \includegraphics[width=0.8\textwidth]{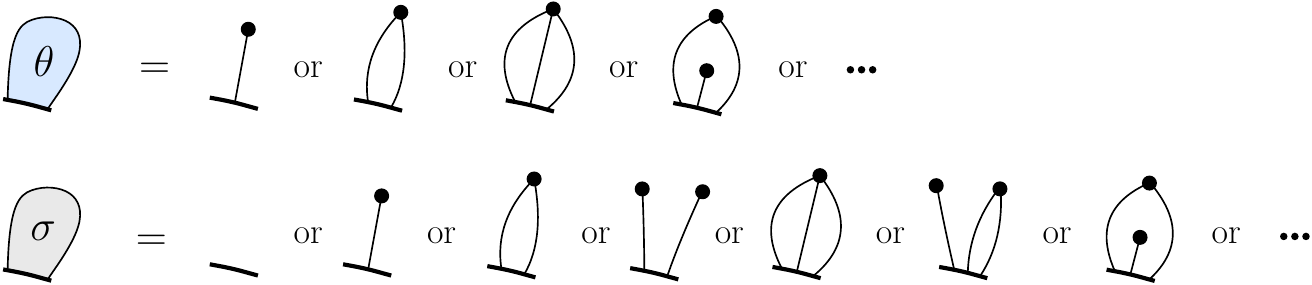}
    \caption{Possibilities for the non-crossing permutations $\sigma^{(\ell)}_{h,i}$ (possibly empty), and $\theta^{(\ell)}_{h,i}$, with up to 3 edges. Letting $n$ be the number of edges of $\theta^{(\ell)}_{h,i}$, one has that $\theta^{(\ell)}_{h,i}(n)=1$.}
    \label{f42}
\end{figure}

\begin{figure}[!h]
    \centering    
    \includegraphics[width=0.65\textwidth]{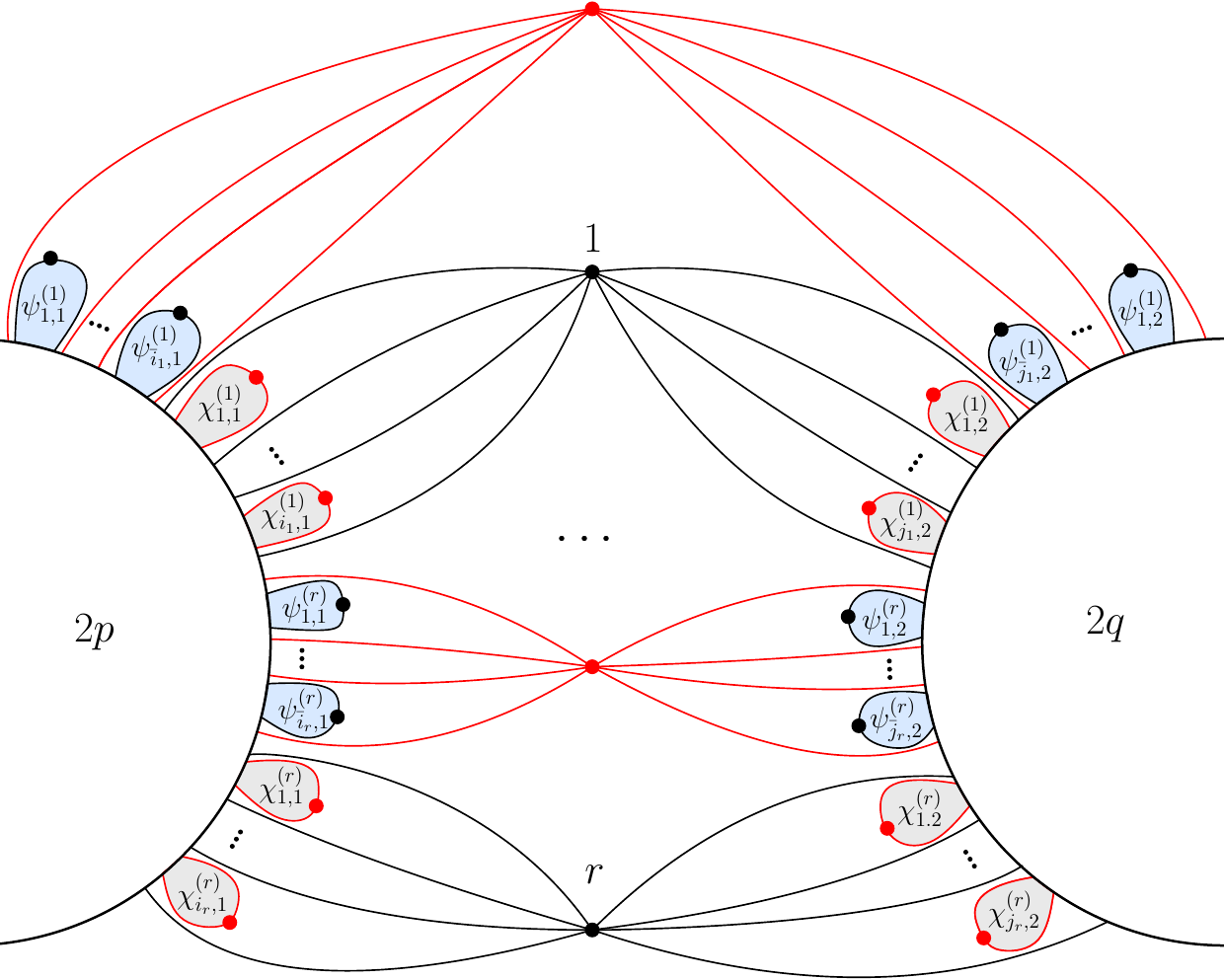}
    \caption{Representation of  $\pi\times Kr(\pi)\in S_{NC}(\operatorname{Int}(\gamma_{p, q}))$, where $\operatorname{Int}(\gamma_{p, q})=(1,\bar{1},\dots,p,\bar{p})(p+1,\overline{p+1},\dots,p+q,\overline{p+q})$. Here, $\psi_{h,s}^{(\ell)}=\theta_{h,s}^{(\ell)}\times Kr^{-1}(\theta_{h,s}^{(\ell)})$ and $\chi_{h,s}^{(\ell)}=\sigma_{h,s}^{(\ell)}\times Kr^{+1}(\sigma_{h,s}^{(\ell)})$}.
 \label{f5} 
\end{figure}

\begin{figure}[!h]
    \centering    
            \includegraphics[scale=0.65]{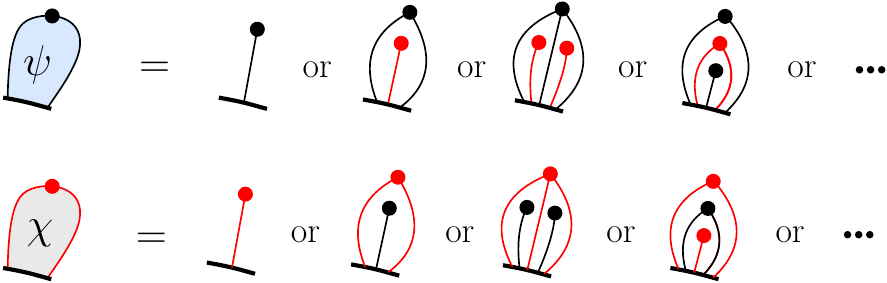}
    \caption{ Due to  Lemma~\ref{lem:Kreweras-prop} (which is equivalent to the  map representation of  $\pi\times Kr(\pi)\in S_{NC}(\operatorname{Int}(\gamma_{p, q}))$ being planar and with only square faces), the contributions $\psi^{(\ell)}_{h,s}$ or $\chi^{(\ell)}_{h,s}$ are respectively of the form $\theta_{h,s}^{(\ell)}\times Kr^{-1}(\theta_{h,s}^{(\ell)})$ and $\sigma_{h,s}^{(\ell)}\times Kr^{+1}(\sigma_{h,s}^{(\ell)})$. Here, $\sigma_{h,s}^{(\ell)}$ is an arbitrary (possibly empty) non-crossing permutation with $p^{(\ell)}_h$ edges if $s=1$ and $q^{(\ell)}_h$ edges if $s=2$, and $\theta_{h,s}^{(\ell)}$ is a non-empty non-crossing permutation satisfying the property that $\theta_{h,s}^{(\ell)}(n+1)=1$, where we have denoted the number of edges of $\theta_{h,s}^{(\ell)}$ by $n+1$. 
    }
 \label{f52} 
\end{figure}

\begin{figure}[!h]
    \centering    
            \includegraphics[scale=0.65]{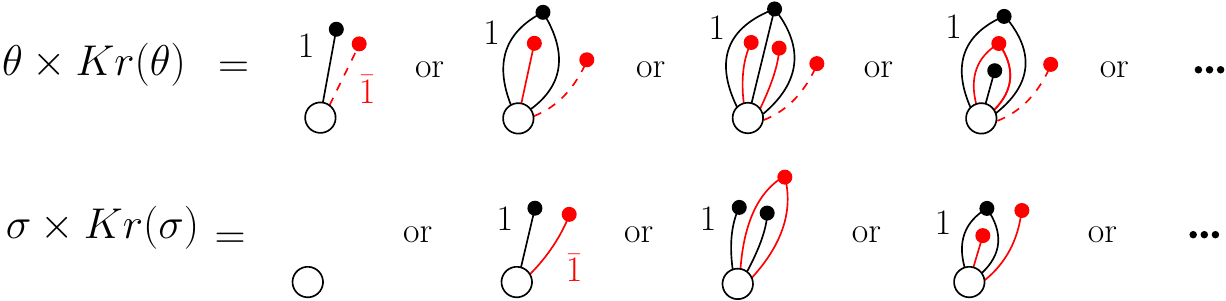}
    \caption{In this figure, we represent $\theta\times Kr(\theta)$ and $\sigma\times Kr(\sigma)$, where $\theta, \sigma\in S_{NC}(n+1)$ and  $\theta$ is such that $\theta(n+1)=1$, for $n=0,1,2$.  In order to recover $\psi=\theta\times Kr^{-1}(\theta)$ in Fig.~\ref{f52}, the red dashed edge (labeled $\overline{n+1}$) must be removed in the first line, and in order to recover $\chi=\sigma\times Kr^{+1}(\sigma)$ in Fig.~\ref{f52}, a red edge labeled $\bar 0$) must be added to the left of the white disc, connected to the red disc which is also connected to the edge
    labeled {$\overline{n+1}$}.  
    }
 \label{f53} 
\end{figure}
\begin{proof}
Since \( \pi \vee \gamma_{p,q} = 1_{p+q} \), it follows that \( \pi \) contains at least one through-cycle. Let \( r \ge 1 \) be the number of through-cycles of \( \pi \in S_{NC}(p,q) \).  For a permutation \( \pi \in S_{NC}(p,q) \) with exactly \( r \) through-cycles, its cycle decomposition thereby consists of precisely \( r \) through-cycles and \( \#\pi - r \) non-through-cycles. 
As illustrated in Fig.~\ref{f4}, we label from 1 to $r$ the through-cycles of $\pi $ according to the minimal label of edges that are connected to them: the through-cycle of $\pi$ that contains the through-edge of minimal label on the first $\gamma_{p,q}$ cycle is labeled as $\ell=1$, the through-cycle that contains the through-edge of minimal label among the edges not connected to the first through-cycle is labeled as $\ell=2$, and so on.  For each $\ell\in\{1, \ldots, r\}$, we let $i_{\ell}+1\ge 1$ and $j_\ell+1\ge 1$ be the numbers of edges that connect the through-cycle number $\ell$ to the first and to the second cycles of $\gamma_{p,q}$, respectively. One has that
\begin{equation}
\label{eq:mn}
\sum_{\ell=1}^r (i_{\ell}+1) = m, \qquad 
\sum_{\ell=1}^r (j_{\ell}+1) = n,
\end{equation}
where  $m\ge r$ and $n\ge r$ are the total numbers of edges that connect the through-cycles to the first and second cycles of \( \gamma_{p,q} \), respectively. 
Once \( r, m, n, \{i_{\ell}\}_{\ell=1}^r, \{j_{\ell}\}_{\ell=1}^r \) are fixed, removing the through-cycles from \( \pi \) leaves a certain number of non-crossing permutations on the first and second cycles of \( \gamma_{p,q}\). 

Some of these non-crossing permutations are located between two consecutive edges connected to the same through-cycle. Going around the first cycle of $\gamma_{p,q}$ clockwise, we consider the first edge encountered and that is connected to the through-cycle labeled 1 (it is not necessarily labeled 1). Between this edge and the second edge connected to the first through-cycle, if the latter exists, there is a (possibly empty) non-crossing permutation $\sigma^{(1)}_{1,1}$. Between this second edge and the next edge encountered clockwise connected to the first through-cycle (if it exists), there is a (possibly empty) non-crossing permutation $\sigma^{(1)}_{2,1}$, and so on. Pursuing this way, we denote (when they exist) by 
\[
\left\{\{\sigma^{(\ell)}_{h,1}\}_{h=1}^{i_{\ell}}\right\}_{\ell=1}^r, 
\qquad 
\left\{\{\sigma^{(\ell)}_{h,2}\}_{h=1}^{j_{\ell}}\right\}_{\ell=1}^r,
\]  the non-crossing permutations located between consecutive through-edges  of the  first and second cycles of \( \gamma_{p,q} \), respectively, 
and that connect to the \emph{same} through-cycle, see Fig. \ref{f42}. On the second cycle of $\gamma_{p,q}$, we use the same convention for the labeling of the $\sigma^{(\ell)}_{h,2}$, but going counterclockwise around this cycle instead of clockwise.\footnote{That is,  going around the second cycle of $\gamma_{p,q}$ counterclockwise, $\sigma^{(1)}_{1,2}$ is located between the first edge  connected to the through-cycle labeled 1 and the second edge connected to the first through-cycle; $\sigma^{(1)}_{2,2}$ is located between this last edge and the third edge  connected to the through-cycle labeled 1, and so on.  } We also denote their cardinalities by
\[
\left\{\{p^{(\ell)}_h\}_{h=1}^{i_{\ell}}\right\}_{\ell=1}^r\geq 0, 
\qquad 
\left\{\{q^{(\ell)}_h\}_{h=1}^{j_{\ell}}\right\}_{\ell=1}^r\geq 0.
\] 
Note that if a single edge connects the through-cycle number $\ell$ to the first cycle of $\gamma_{p,q}$, then there is no $\sigma^{(\ell)}_{h,1}$.

The non-crossing permutations listed above do not include the ones located between consecutive through-edges belonging to different through-cycles. Because we need to associate weights to the cycles of the Kreweras complement $Kr(\pi)$ of $\pi$, we need in particular to tract the number of edges connected to the through-cycles of $Kr(\pi)$. We therefore partition the remaining non-crossing permutations as follows. 

Turning around the first cycle of $\gamma_{p,q}$ clockwise, we  focus on the (possibly empty) non-crossing permutation located between the last edge connected to the through-cycle labeled $\ell-1$, and the first edge connected to the through-cycle labeled $\ell$ (when $\ell=1$, $\ell -1$ is replaced by $r$). This non-crossing permutation is itself partitioned into a certain number $\bar i_\ell\ge 0$ of non-crossing permutations that have the additional property that the first edge (labeled $e\in \{1, \ldots, p+q\}$) and the  last edge (labeled $e'$) encountered clockwise when going around the first cycle of $\gamma_{p,q}$ are connected to the same cycle of $\pi$. Equivalently, $\pi$ has the property that $\pi(e')=e$. We denote by $\{\theta^{(\ell)}_{h,1}\}_{h=1}^{\bar{i}_{\ell}}$ these $\bar i_\ell$ non-crossing permutations, where $h$ grows when going clockwise around the first cycle of $\gamma_{p,q}$.

Going counterclockwise around the second cycle of $\gamma_{p,q}$, we define in the same way the numbers $\bar j_\ell\ge 0$ and the non-crossing permutations $\{\theta^{(\ell)}_{h,2}\}_{h=1}^{\bar{j}_{\ell}}$, satisfying the same property, see Fig. \ref{f42}. 
The number of elements on which the non-crossing permutations  
\[
\left\{\{\theta^{(\ell)}_{h,1}\}_{h=1}^{\bar{i}_{\ell}}\right\}_{\ell=1}^r, 
\qquad 
\left\{\{\theta^{(\ell)}_{h,2}\}_{h=1}^{\bar{j}_{\ell}}\right\}_{\ell=1}^r,
\] act (the number of edges that they involve) are denoted by
\[
\left\{\{\bar{p}^{(\ell)}_h+1\}_{h=1}^{\bar{i}_\ell}\right\}_{\ell=1}^r \geq 1, 
\qquad
\left\{\{\bar{q}^{(\ell)}_h+1\}_{h=1}^{\bar{j}_\ell}\right\}_{\ell=1}^r \geq 1, 
\]
respectively. 
 By construction,\footnote{With the convention that if $i_\ell=0$ of $\bar i_\ell=0$, the corresponding sum is not considered in the equation below. }
\begin{equation} \label{eq:q-n}
\sum_{\ell=1}^r \left( \sum_{h=1}^{i_{\ell}} p^{(\ell)}_h 
      + \sum_{h=1}^{\bar{i}_{\ell}} (\bar{p}^{(\ell)}_h+1) \right)
    = p - m, \quad 
\sum_{\ell=1}^r \left( \sum_{h=1}^{j_{\ell}} q^{(\ell)}_h
      + \sum_{h=1}^{\bar{j}_{\ell}} (\bar{q}^{(\ell)}_h+1) \right)
    = q - n.
\end{equation}

Fixing all previously defined parameters and choosing the initial element of the first cycle of \( \gamma_{p,q} \) (the edge labeled 1) and the initial element of the second cycle of \( \gamma_{p,q} \) (the edge of labeled $p+1$) completely determines $\pi$ (that is, it is a one-to-one correspondence).

Furthermore, due  to  Lemma~\ref{lem:Kreweras-prop} (which is equivalent graphically to the  map representation of $\bigl(\operatorname{Int}(\gamma_{p,q}), \pi\times Kr(\pi)\bigr)$ being planar and with only square faces), and to the fact that the non-crossing permutations $\theta^{(\ell)}_{h,1}$ and the $\theta^{(\ell)}_{h,2}$ have the additional property stated above,  between the through-cycle of $\pi$ labeled $\ell$ and that labeled $\ell+1$, there is precisely one through-cycle of $Kr(\pi)$, and it has precisely $\bar i_\ell + \bar j_\ell + 2$ edges attached,  see Figures \ref{f5}, \ref{f52} and  \ref{f53}. In addition to this, focusing on two through-edges of $\pi$ connected to the same through-cycle and which bound a non-crossing permutation $\sigma^{(\ell)}_{h, s}$, $s\in \{1,2\}$, then  due again to Lemma~\ref{lem:Kreweras-prop}, the contribution to the Kreweras complement of $\pi$ that is located between these two consecutive through-edges is  $Kr^{+1}(\sigma_{h,s}^{(\ell)})$ (this notation has been introduced in Prop.~\ref{p1}), as illustrated in Fig. \ref{f53}. In the same way, the contributions to the Kreweras complement of $\pi$ that are located between the edges of  $\theta^{(\ell)}_{h, s}$, $s\in \{1,2\}$, are given by $Kr^{-1}(\theta_{h,s}^{(\ell)})$.
Therefore:
\[
\begin{aligned}
\kappa^{a}_{\pi}\kappa^{b}_{K r(\pi)}
&=\prod_{\ell=1}^r\Bigg[
\kappa^{a}_{i_{\ell}+j_{\ell}+2}\,\kappa^{b}_{\bar{i}_{\ell}+\bar{j}_{\ell}+2}
\prod_{h=1}^{i_{\ell}} \kappa^{a}_{\sigma_{h,1}^{(\ell)}}  \kappa^{b}_{Kr^{+1}(\sigma_{h,1}^{(\ell)})}
\prod_{h=1}^{j_{\ell}} \kappa^{a}_{\sigma_{h,2}^{(\ell)}}  \kappa^{b}_{Kr^{+1}(\sigma_{h,2}^{(\ell)})} 
\prod_{h=1}^{\bar{i}_{\ell}}\kappa^{a}_{\theta_{h,1}^{(\ell)}}  \kappa^{b}_{Kr^{-1}(\theta_{h,1}^{(\ell)})}
\prod_{h=1}^{\bar{j}_{\ell}} \kappa^{a}_{\theta_{h,2}^{(\ell)}}  \kappa^{b}_{Kr^{-1}(\theta_{h,2}^{(\ell)})}
\Bigg].
\end{aligned}
\]

 Thus, applying  the statements $(c), (e)$ and $(f)$ of Proposition \ref{p1}, the sum over all independent non-crossing permutations satisfying the previously imposed conditions yields the contribution
\begin{equation}
\label{eq:unlabeled-contribution}
\prod_{\ell=1}^r\Bigg[
\kappa^{a}_{i_{\ell}+j_{\ell}+2}\,\kappa^{b}_{\bar{i}_{\ell}+\bar{j}_{\ell}+2}
\prod_{h=1}^{i_{\ell}} \kappa_{b(ab)^{p_h^{(\ell)}}}
\prod_{h=1}^{j_{\ell}} \kappa_{b(ab)^{q_h^{(\ell)}}}
\prod_{h=1}^{\bar{i}_{\ell}}\kappa_{a(ba)^{\bar{p}_h^{(\ell)}}}
\prod_{h=1}^{\bar{j}_{\ell}}\kappa_{a(ba)^{\bar{q}_h^{(\ell)}}}
\Bigg].
\end{equation}

Finally, observe that since by convention the through-cycle of $\pi$ that contains the through-edge of minimal label on the first $\gamma_{p,q}$ cycle is labeled as $\ell=1$, the edge labeled 1 must be either a through-edge connected to the through-cycle $\ell=1$, or it is in one of the $\sigma_{h,1}^{(1)}$, or it is in one of the $\theta_{h,1}^{(1)}$. Therefore, the number of ways of choosing the position of the edge labeled 1 is:
\[
i_1 + 1 + \sum_{h=1}^{i_{1}} p^{(1)}_{h}+ \sum_{h=1}^{\bar{i}_{1}}( \bar{p}^{(1)}_{h}+1). 
\]
On the other hand, there are \(q\) possible  distinct choices for the edge labeled $p+1$.
 Consequently, there are 
\(q ~ (i_1 + 1 + \sum_{h=1}^{i_{1}} p^{(1)}_h + \bar{i}_{1} + \sum_{h=1}^{\bar{i}_{1}} \bar{p}^{(1)}_{h})\) 
distinct possibilities for \( \pi \) that all give the same  contribution \eqref{eq:unlabeled-contribution}. Therefore, we obtain

\begin{align*}
{\kappa}^{ab,1}_{p,q}
&=
\sum_{r \ge 1}\sum_{m,n \ge r}
\sum_{\substack{
   \{i_\ell, j_\ell\}_{\ell=1}^r \ge 0 \\
   \sum_{\ell=1}^r i_\ell = m - r \\
   \sum_{\ell=1}^r j_\ell = n - r
}}
\sum_{\substack{
  { \{\bar{i}_\ell, \bar{j}_\ell\}_{\ell=1}^r \ge 0 }\\{
  \left\{\{p^{(\ell)}_h,q^{(\ell)}_h\}_{h=1}^{i_{\ell},j_{\ell}}\right\}_{\ell=1}^r \ge 0 }\\{
    \left\{\{\bar{p}^{(\ell)}_h,\bar{q}^{(\ell)}_h\}_{h=1}^{\bar{i}_{\ell},\bar{j}_{\ell}}\right\}_{\ell=1}^r \ge 0}\\[4pt]
    {\sum_{\ell=1}^r \sum_{h=1}^{i_\ell} p^{(\ell)}_h 
      + \bar{i}_\ell 
      + \sum_{h=1}^{\bar{i}_\ell} \bar{p}^{(\ell)}_h= p - m} 
      \\[3pt]
   {\sum_{\ell=1}^r \sum_{h=1}^{j_\ell} q^{(\ell)}_h
      + \bar{j}_\ell
      + \sum_{h=1}^{\bar{j}_\ell} \bar{q}^{(\ell)}_h = q - n}
      }
      }
q\left(i_1 + 1 + \sum_{h=1}^{i_{1}} p^{(1)}_h + \bar{i}_{1} + \sum_{h=1}^{\bar{i}_{1}} \bar{p}^{(1)}_{h} \right)
\\[4pt]
&\qquad\qquad\qquad
\cdot \prod_{\ell=1}^r
\Bigg[
\kappa^{a}_{i_{\ell} + j_{\ell} + 2}\,
\kappa^{b}_{\bar{i}_{\ell} + \bar{j}_{\ell} + 2} 
\cdot\prod_{h=1}^{i_{\ell}} \kappa_{b(ab)^{\,p_h^{(\ell)}}}
\prod_{h=1}^{j_{\ell}} \kappa_{b(ab)^{\,q_h^{(\ell)}}}
\prod_{h=1}^{\bar{i}_{\ell}} \kappa_{a(ba)^{\,\bar{p}_h^{(\ell)}}}
\prod_{h=1}^{\bar{j}_{\ell}} \kappa_{a(ba)^{\,\bar{q}_h^{(\ell)}}}
\Bigg].
\end{align*}
 Hence, using \eqref{eq:mn} and \eqref{eq:q-n} to  express $p,q$ in terms of the other parameters:
\[
\begin{aligned}
\sum_{p,q\ge 1}\kappa^{ab,1}_{p,q}\, x^{p-1}y^{q-1}
&=
\sum_{r \ge 1} \sum_{\substack{
   \{i_\ell, j_\ell\}_{\ell=1}^r \ge 0 \\
   \{\bar{i}_\ell, \bar{j}_\ell\}_{\ell=1}^r \ge 0 
}}
\sum_{\substack{
    \left\{\{p^{(\ell)}_h,q^{(\ell)}_h\}_{h=1}^{i_{\ell},j_{\ell}}\right\}_{\ell=1}^r \ge 0 \\
    \left\{\{\bar{p}^{(\ell)}_h,\bar{q}^{(\ell)}_h\}_{h=1}^{\bar{i}_{\ell},\bar{j}_{\ell}}\right\}_{\ell=1}^r \ge 0 
  }}
\left(i_1 + 1 + \sum_{h=1}^{i_{1}} p^{(1)}_h + \bar{i}_{1} + \sum_{h=1}^{\bar{i}_{1}} \bar{p}^{(1)}_{h} \right)\\[4pt] &\qquad\qquad\qquad
\cdot 
\left( \sum_{\ell=1}^r  \Bigl( j_{\ell}+1+\sum_{h=1}^{j_{\ell}} q^{(\ell)}_h
      + \bar{j}_\ell
      + \sum_{h=1}^{\bar{j}_{\ell}} \bar{q}^{(\ell)}_h \Bigr)\right)\\[4pt]
&\qquad\qquad\qquad
\cdot  \prod_{\ell=1}^r
\Bigg[
\kappa^{a}_{i_{\ell} + j_{\ell} + 2}\,
\kappa^{b}_{\bar{i}_{\ell} + \bar{j}_{\ell} + 2} 
\cdot\prod_{h=1}^{i_{\ell}} \kappa_{b(ab)^{\,p_h^{(\ell)}}}
\prod_{h=1}^{j_{\ell}} \kappa_{b(ab)^{\,q_h^{(\ell)}}}
\prod_{h=1}^{\bar{i}_{\ell}} \kappa_{a(ba)^{\,\bar{p}_h^{(\ell)}}}
\prod_{h=1}^{\bar{j}_{\ell}} \kappa_{a(ba)^{\,\bar{q}_h^{(\ell)}}}
\Bigg]\\[4pt]
&\qquad\qquad\qquad
\cdot x^{\sum_{\ell=1}^r \left( i_{\ell}+1+\sum_{h=1}^{i_{\ell}}p^{(\ell)}_h
      + \bar{i}_\ell
      + \sum_{h=1}^{\bar{i}_{\ell}} \bar{p}^{(\ell)}_h \right)-1}y^{\sum_{\ell=1}^r \left( j_{\ell}+1+\sum_{h=1}^{j_{\ell}} q^{(\ell)}_h
      + \bar{j}_\ell
      + \sum_{h=1}^{\bar{j}_{\ell}} \bar{q}^{(\ell)}_h \right)-1}\\[4pt]
&
=\sum_{r \ge 1} \frac{\partial}{\partial y} \Bigg[  \sum_{\substack{
   \{i_\ell, j_\ell\}_{\ell=1}^r \ge 0 \\
   \{\bar{i}_\ell, \bar{j}_\ell\}_{\ell=1}^r \ge 0 
}}
\sum_{\substack{
    \left\{\{p^{(\ell)}_h,q^{(\ell)}_h\}_{h=1}^{i_{\ell},j_{\ell}}\right\}_{\ell=1}^r \ge 0 \\
    \left\{\{\bar{p}^{(\ell)}_h,\bar{q}^{(\ell)}_h\}_{h=1}^{\bar{i}_{\ell},\bar{j}_{\ell}}\right\}_{\ell=1}^r \ge 0
  }}
\left(i_1 + 1 + \sum_{h=1}^{i_{1}} p^{(1)}_h + \bar{i}_{1} + \sum_{h=1}^{\bar{i}_{1}} \bar{p}^{(1)}_{h}  \right)\\[4pt]
&\qquad\qquad\qquad
 \cdot \prod_{\ell=1}^r
\Bigg[
\kappa^{a}_{i_{\ell} + j_{\ell} + 2}\,
\kappa^{b}_{\bar{i}_{\ell} + \bar{j}_{\ell} + 2} 
\cdot\prod_{h=1}^{i_{\ell}} \kappa_{b(ab)^{\,p_h^{(\ell)}}}
\prod_{h=1}^{j_{\ell}} \kappa_{b(ab)^{\,q_h^{(\ell)}}}
\prod_{h=1}^{\bar{i}_{\ell}} \kappa_{a(ba)^{\,\bar{p}_h^{(\ell)}}}
\prod_{h=1}^{\bar{j}_{\ell}} \kappa_{a(ba)^{\,\bar{q}_h^{(\ell)}}}
\Bigg]\\[4pt]
&\qquad\qquad\qquad
\cdot x^{\sum_{\ell=1}^r \left( i_{\ell}+1+\sum_{h=1}^{i_{\ell}}p^{(\ell)}_h
      + \bar{i}_\ell
      + \sum_{h=1}^{\bar{i}_{\ell}} \bar{p}^{(\ell)}_h \right)-1}y^{\sum_{\ell=1}^r \left( j_{\ell}+1+\sum_{h=1}^{j_{\ell}} q^{(\ell)}_h
      + \bar{j}_\ell
      + \sum_{h=1}^{\bar{j}_{\ell}} \bar{q}^{(\ell)}_h \right)} \Bigg].
\end{aligned}
\] Now, note that the expression inside the derivative is equivalent to

\[
\begin{aligned}
&
\sum_{\substack{
   \{i_\ell, j_\ell\}_{\ell=1}^r \ge 0 \\
   \{\bar{i}_\ell, \bar{j}_\ell\}_{\ell=1}^r \ge 0
}}
\Bigg[
\prod_{\ell=1}^r
\kappa^{a}_{i_\ell + j_\ell + 2}\,
\kappa^{b}_{\bar{i}_\ell + \bar{j}_\ell + 2}
\\[6pt]
&\quad \cdot
\sum_{\substack{
   \{p^{(1)}_h\}_{h=1}^{i_1} \ge 0 \\
   \{\bar{p}^{1}_h\}_{h=1}^{\bar{i}_1} \ge 0
}}\Bigg(
\left(
i_1 + 1 + \sum_{h=1}^{i_{1}} p^{(1)}_h + \bar{i}_{1} + \sum_{h=1}^{\bar{i}_{1}} \bar{p}^{(1)}_{h} 
\right)
\prod_{h=1}^{i_1}
\kappa_{b(ab)^{p_h^{(1)}}}
\prod_{h=1}^{\bar{i}_1}
\kappa_{a(ba)^{\bar{p}_h^{(1)}}}
x^{
 i_1 + 1
 + \sum_{h=1}^{i_1} p^{(1)}_h
 + \bar{i}_1
 + \sum_{h=1}^{\bar{i}_1} \bar{p}^{(1)}_h
-1}\Bigg)
\\[6pt]
&\quad \cdot
\prod_{\ell=2}^r
\Bigg(
  \prod_{h=1}^{i_\ell}
  \sum_{p_h^{(\ell)} \ge 0}
  \kappa_{b(ab)^{p_h^{(\ell)}}}
  x^{p_h^{(\ell)}}
\Bigg)
x^{\,i_\ell+1}
\quad \cdot
\prod_{\ell=1}^r
\Bigg(
  \prod_{h=1}^{j_\ell}
  \sum_{q_h^{(\ell)} \ge 0}
  \kappa_{b(ab)^{q_h^{(\ell)}}}
  y^{q_h^{(\ell)}}
\Bigg)
y^{\,j_\ell+1}
\\[6pt]
&\quad \cdot
\prod_{\ell=2}^r
\Bigg(
  \prod_{h=1}^{\bar{i}_\ell}
  \sum_{\bar{p}_h^{(\ell)} \ge 0}
  \kappa_{a(ba)^{\bar{p}_h^{(\ell)}}}
  x^{\bar{p}_h^{(\ell)}}
\Bigg)
x^{\,\bar{i}_\ell}
\quad \cdot
\prod_{\ell=1}^r
\Bigg(
  \prod_{h=1}^{\bar{j}_\ell}
  \sum_{\bar{q}_h^{(\ell)} \ge 0}
  \kappa_{a(ba)^{\bar{q}_h^{(\ell)}}}
  y^{\bar{q}_h^{(\ell)}}
\Bigg)
y^{\,\bar{j}_\ell}
\Bigg].
\end{aligned}
\] By introducing the generating functions
\[
C_1(z)=\sum_{n\ge 0} \kappa_{b(ab)^n} \, z^n,
\qquad
C_2(z)=\sum_{n\ge 0} \kappa_{a(ba)^n}\, z^n,
\]
the preceding expression can be rewritten in the following form.

\[
\begin{aligned}
&
\sum_{\substack{
   \{i_\ell, j_\ell\}_{\ell=1}^r \ge 0 \\
   \{\bar{i}_\ell, \bar{j}_\ell\}_{\ell=1}^r \ge 0
}}
\Bigg[
\prod_{\ell=1}^r
\kappa^{a}_{i_\ell + j_\ell + 2}\,
\kappa^{b}_{\bar{i}_\ell + \bar{j}_\ell + 2}
\\[6pt]
&\quad \cdot
\sum_{\substack{
   \{p^{(1)}_h\}_{h=1}^{i_1} \ge 0 \\
   \{\bar{p}^{(1)}_h\}_{h=1}^{\bar{i}_1} \ge 0
}}\Bigg(
\left(
i_1 + 1
+ \sum_{h=1}^{i_1} p^{(1)}_h
+ \bar{i}_1
+ \sum_{h=1}^{\bar{i}_1} \bar{p}^{(1)}_h
\right)
\prod_{h=1}^{i_1}
\kappa_{b(ab)^{p_h^{(1)}}}
\prod_{h=1}^{\bar{i}_1}
\kappa_{a(ba)^{\bar{p}_h^{(1)}}}
x^{
 i_1 + 1
 + \sum_{h=1}^{i_1} p^{(1)}_h
 + \bar{i}_1
 + \sum_{h=1}^{\bar{i}_1} \bar{p}^{(1)}_h
-1}\Bigg)
\\[6pt]
&\quad \cdot
\prod_{\ell=2}^r
\left(
   \left( x\,C_1(x) \right)^{i_{\ell+1}}
   \frac{1}{C_1(x)}
\right)
 \cdot
\prod_{\ell=1}^r
\left(
   \left( y\,C_1(y) \right)^{j_{\ell+1}}
   \frac{1}{C_1(y)}
\right) \cdot
\prod_{\ell=2}^r
\left(
   \left( x\,C_2(x) \right)^{\bar{i}_{\ell}}
\right)
\cdot
\prod_{\ell=1}^r
\left(
   \left( y\,C_2(y) \right)^{\bar{j}_{\ell}}
\right)
\Bigg].
\end{aligned}
\]  On the other hand,

\[
\begin{aligned}
&\sum_{\substack{
   \{ p^{(1)}_h \}_{h=1}^{i_1} \ge 0 \\
   \{ \bar{p}^{(1)}_h \}_{h=1}^{\bar{i}_1} \ge 0 
}} \Bigg(
\left(
i_1 + 1 
+ \sum_{h=1}^{i_1} p^{(1)}_h
+ \bar{i}_1
+ \sum_{h=1}^{\bar{i}_1} \bar{p}^{(1)}_h
\right)
\prod_{h=1}^{i_1} 
\kappa_{b(ab)^{p_h^{(1)}}}
\prod_{h=1}^{\bar{i}_1} 
\kappa_{a(ba)^{\bar{p}_h^{(1)}}}
x^{
 i_1 + 1 
 + \sum_{h=1}^{i_1} p^{(1)}_h
 + \bar{i}_1
 + \sum_{h=1}^{\bar{i}_1} \bar{p}^{(1)}_h
-1}
 \Bigg)
\\[6pt]
&\qquad
=\frac{\partial}{\partial x}
\Bigg[
\Bigg(
\prod_{h=1}^{i_1}
\sum_{\substack{
    p^{(1)}_h \ge 0
}}
\kappa_{b(ab)^{p_h^{(1)}}} \,
x^{p_h^{(1)}}
\Bigg)
\Bigg(
\prod_{h=1}^{\bar{i}_1}
\sum_{\substack{
    \bar{p}^{(1)}_h \ge 0
}}
\kappa_{a(ba)^{\bar{p}_h^{(1)}}} \,
x^{\bar{p}_h^{(1)}}
\Bigg) x^{i_1+\bar{i}_1+1}
\Bigg]
\\[8pt]
&\qquad
= \frac{\partial}{\partial x}
\Bigg[
\frac{ \big( x\, C_1(x)\big)^{\,i_1+1} \, \big( x\, C_2(x)\big)^{\,\bar{i}_1} }
       { C_1(x) }
\Bigg].
\end{aligned}
\] Thus, after substituting the equality above, the quantity inside the partial derivative with respect to $y$ takes the form

\[
\begin{aligned}
&
\frac{\partial}{\partial x}
\Bigg[
\frac{1}{C_1(x)\,C_1(y)}
\sum_{\substack{i_1,j_1 \ge 0}}
\kappa^{a}_{i_1 + j_1 + 2}\,
\left( x\,C_1(x) \right)^{i_1+1}
\left( y\,C_1(y) \right)^{j_1+1}
\cdot
\sum_{\substack{\bar{i}_1,\bar{j}_1 \ge 0}}
\kappa^{b}_{\bar{i}_1 + \bar{j}_1 + 2}\,
\left( x\,C_2(x) \right)^{\bar{i}_1}
\left( y\,C_2(y) \right)^{\bar{j}_1}
\Bigg]
\\[6pt]
&\quad\cdot
\Bigg[
\frac{1}{(C_1(x)\,C_1(y))^{r-1}}
\prod_{\ell=2}^r
\sum_{\substack{i_{\ell},j_{\ell}\ge 0}}
\kappa^{a}_{i_{\ell}+j_{\ell}+2}
\left( x\,C_1(x) \right)^{i_{\ell}+1}
\left( y\,C_1(y) \right)^{j_{\ell}+1}
\Bigg]
\Bigg[
\prod_{\ell=2}^r
\sum_{\substack{\bar{i}_{\ell},\bar{j}_{\ell}\ge 0}}
\kappa^{b}_{\bar{i}_{\ell}+\bar{j}_{\ell}+2}
\left( x\,C_2(x) \right)^{\bar{i}_{\ell}}
\left( y\,C_2(y) \right)^{\bar{j}_{\ell}}
\Bigg].
\end{aligned}
\] Now, considering 
\[ \widehat C_d(x,y)=\sum_{i,j\ge 0}\kappa_{i+j+2}^{\,d}\,x^{\,i+1}y^{\,j+1},
\qquad d\in\{a,b\},
\] we have that

\[
\begin{aligned}
R_{ab,1}=\sum_{p,q\ge 1}\kappa^{ab,1}_{p,q}\, x^{p-1}y^{q-1}
&= \sum_{r \ge 1}\Bigg[ \frac{\partial}{\partial y} \Bigg[
\frac{\partial}{\partial x}
\Bigg[
\frac{\widehat{C}_a(x\,C_1(x), y\,C_1(y))\,\widehat{C}_b(x\,C_2(x), y\,C_2(y))}{x\,y\,C_1(x)\,C_1(y)\,C_2(x)\,C_2(y)}
\Bigg]
\\[6pt]
&\hspace{2cm}\cdot
\Bigg(\frac{\widehat{C}_a(x\,C_1(x), y\,C_1(y))}{C_1(x)\,C_1(y)}\Bigg)^{r-1}
\cdot
\Bigg(\frac{\widehat{C}_b(x\,C_2(x), y\,C_2(y))}{x\,y\,C_2(x)\,C_2(y)}\Bigg)^{r-1}\Bigg] \Bigg].
\end{aligned}
\] Finally, 
\[
\begin{aligned}
R_{ab,1}(x,y)
&= \frac{\partial}{\partial y}
   \Bigg[
      \frac{\partial}{\partial x}
      \Biggl[
         {\frac{
            \widehat{C}_a\bigl(x\,C_1(x),\, y\,C_1(y)\bigr)\,
            \widehat{C}_b\bigl(x\,C_2(x),\, y\,C_2(y)\bigr)
         }{
            x\,y\,
            C_1(x)\,C_1(y)\,
            C_2(x)\,C_2(y)
            }
            }
      \Biggr] \cdot
      {\frac{
         1
      }{
         1-
         \dfrac{
            \widehat{C}_a\bigl(x\,C_1(x),\, y\,C_1(y)\bigr)\,
            \widehat{C}_b\bigl(x\,C_2(x),\, y\,C_2(y)\bigr)
         }{
            x\,y\,
            C_1(x)\,C_1(y)\,
            C_2(x)\,C_2(y)
         }
    }}
   \Bigg], 
\end{aligned}
\] 
which we rewrite as
\[
R_{ab,1}(x,y)
= - \frac{\partial^2}{\partial y \partial x}
         \log\Biggl( 1 - \frac{
            \widehat{C}_a\bigl(x\,C_1(x),\, y\,C_1(y)\bigr)\,
            \widehat{C}_b\bigl(x\,C_2(x),\, y\,C_2(y)\bigr)
         }{
            x\,y\,
            C_1(x)\,C_1(y)\,
            C_2(x)\,C_2(y)
            }
            \Biggr). 
\] Therefore, by \eqref{ecr} the proof is complete. \end{proof}

\begin{lem}\label{l2}
Let $a_1$ and $b_1$ be a pair of first-order free random variables and
let $a_2$ and $b_2$ be another pair of first-order free
random variables. Let $p,q\geq 1$ be integers. Define 
\begin{equation}
\label{eq:combinatorial-factorized-product1}
\begin{aligned}
\kappa_{p,q}^{\vec a, \vec b , 2}
&=
\sum_{\substack{
\pi=\pi_1\times\pi_2\\
\pi_1\in S_{NC}(p),\ \pi_2\in S_{NC}(q)
}}
\ \sum_{\substack{U_1\in\pi_1\\ U_2\in\pi_2}}
\kappa_{|U_1|,|U_2|}^{a_1}
\kappa^{ a_1}_{\pi_1\setminus\{U_1\}}
\kappa^{a_2}_{\pi_2\setminus\{U_2\}}
\kappa^{ b_1}_{Kr(\pi_1)}
\kappa^{b_2}_{Kr(\pi_2)}.
\end{aligned}
\end{equation}
Letting:
\[
R_{\vec a, \vec b, 2}(x,y)
=
\sum_{p,q\geq 1}
\kappa_{p,q}^{\vec a, \vec b , 2}\,x^{p-1}y^{q-1}.
\]
and for $i=1,2$:
\[
C_1^{(i)}(z)
=
\sum_{n\geq 0}
\kappa_{b_i(a_ib_i)^n}\,z^n,
\qquad
C_2^{(i)}(z)
=
\sum_{n\geq 0}
\kappa_{a_i(b_ia_i)^n}\,z^n,
\]
the generating function associated with
\eqref{eq:combinatorial-factorized-product1} admits the expression:
\begin{equation}
\label{eq:functional-factorized-product2}
R_{\vec a, \vec b, 2}(x,y)
=
\big(xC_1^{(1)}(x)\big)'
\big(yC_1^{(2)}(y)\big)'
R_a\big(xC_1^{(1)}(x),yC_1^{(2)}(y)\big).
\end{equation}
If the roles of $a$ and $b$ are exchanged in
\eqref{eq:combinatorial-factorized-product1}, then $C_1^{(i)}$ is replaced
by $C_2^{(i)}$ in \eqref{eq:functional-factorized-product2}.
\end{lem}

\begin{proof}

With these notations we consider $\pi=\pi_1\times \pi_2$,  $U_1\in \pi_1$ and $U_2\in \pi_2$, an we let $n\geq 1$ be the length of $U_1$ and $m \geq 1$ that of $U_2$. Removing these two cycles from $\pi$ leaves \( n \) and \( m \)  independent (possibly empty) non-crossing permutations in the first and second cycles of \( \gamma_{p,q} \).  We denote by $\sigma_{1,1}$ (resp.~$\sigma_{1,2}$) the non-crossing permutation located on the first (resp.~second) cycle of $\gamma_{p,q}$, before the edge of minimal label connected to $U_1$ (resp.~$U_2$), clockwise. We then denote the other $n-1$ (resp.~$m-1$) non-crossing permutations as $\{\sigma_{h,1}\}_{h=2}^{n}$ (resp.~$\{\sigma_{h,2}\}_{h=2}^{m}$) with $h$ growing according to their order of appearance when pursuing  clockwise around the first (resp.~second) cycle of $\gamma_{p,q}$.   
We also let $p_h$ (resp.~$q_h$) be the cardinality of $\sigma_{h,1}$ (resp.~$\sigma_{h,2}$). One has:
\[
\sum_{h=1}^n p_{h} = p -n, \qquad 
\sum_{h=1}^m q_{h} = q-m.
\] 
Under our convention, the edge labeled 1 (resp.~$p+1$) is either connected to the cycle $U_1$ (resp.~$U_2$), or it belongs to $\sigma_{1,1}$ (resp.~$\sigma_{1,2}$), see Fig. \ref{f8}. 
 Fixing all previously defined parameters and choosing the initial element of the first cycle of \( \gamma_{p,q} \) (the edge labeled 1) and the initial element of the second cycle of \( \gamma_{p,q} \) (the edge labeled  $p+1$) completely determines $\pi$ (it is a one-to-one correspondence).

\begin{figure}
    \centering
   \qquad\qquad\qquad \quad\ \ \includegraphics[width=0.85\textwidth]{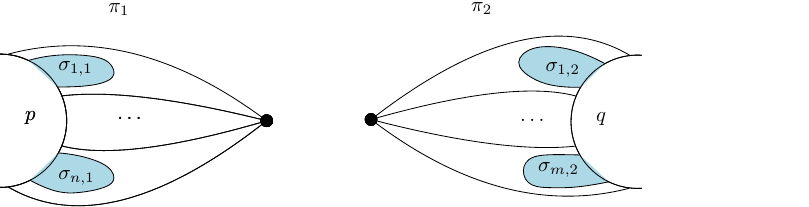}
    \caption{Representation of $\pi_1 \in S_{NC}(p)$ and $\pi_2 \in S_{NC}(q)$ appearing in Lemma \ref{l2}.}
    \label{f8}
\end{figure}

Similarly to the proof of  Lemma~\ref{l1}, the  contribution to the Kreweras complement of $\pi$ that is located just before the edge of minimal label connected to $U_1$ on the first cycle of $\gamma_{p,q}$, clockwise, is identified as $Kr^{+1}(\sigma_{1,1})$, and similarly for the other non-crossing permutations, see Fig. \ref{f53}.  Consequently,
\[
\begin{aligned}
\kappa_{|U_1|,|U_2|}^{a_1}
\kappa^{ a_1}_{\pi_1\setminus\{U_1\}}
\kappa^{ a_2}_{\pi_2\setminus\{U_2\}}
\kappa^{ b_1}_{Kr(\pi_1)}
\kappa^{ b_2}_{Kr(\pi_2)}
&=\kappa_{n,m}^{a_1}
\prod_{h=1}^{n}\kappa^{ a_1}_{\sigma_{h,1}}  \kappa^{ b_1}_{Kr^{+1}(\sigma_{h,1})}
\prod_{h=1}^{m}\kappa^{ a_2}_{\sigma_{h,2}}  \kappa^{ b_2}_{Kr^{+1}(\sigma_{h,2})}
.
\end{aligned}
\] Thus, by Proposition \ref{p1} $(e)$, the sum over all independent non-crossing permutations satisfying the previously imposed conditions yields the contribution
\begin{equation}
\label{eq:second-proof}
\kappa_{n,m}^{a_1}\
\prod_{h=1}^{n} \kappa_{b_1(a_1b_1)^{p_{h}}}
\prod_{h=1}^{m} \kappa_{b_2(a_2b_2)^{q_{h}}}.
\end{equation}
Finally, observe that there are $p_1+1$ 
possible choices for the positions of the edge labeled 1: either it is connected to $U_1$ and this fact entirely fixes $\pi_1$, or it is in $\sigma_{1,1}$, in which case one has to specify which of the $p_1$ edges of $\sigma_{1,1}$ it is. Similarly, there are   $q_1+1$ possible choices for the position of the edge labeled $p+1$. In total, there are therefore
\( (p_1 + 1) (q_1+1)\) distinct possibilities for $(\pi, U_1, U_2)$ that all give the same contribution \eqref{eq:second-proof}. Consequently:
\[
\begin{aligned}
\kappa_{p,q}^{\vec a, \vec b , 2}
&=
\sum_{m,n \ge 1}
\sum_{\substack{
   \{p_h\}_{h=1}^n \ge 0 \\
   \sum_{h=1}^n p_h = p - m \\
 }}
\sum_{\substack{
   \{q_h\}_{h=1}^m \ge 0 \\
   \sum_{h=1}^m q_h = q - n \\
 }}
\left(p_1 + 1\right)\left(q_1 + 1\right)\kappa^{a_1}_{n,m}\
\prod_{h=1}^{n} \kappa_{b_1(a_1b_1)^{p_{h}}}
\prod_{h=1}^{m} \kappa_{b_2(a_2b_2)^{q_{h}}}.
\end{aligned}
\] Hence,

\[
\begin{aligned}
\sum_{p,q\ge 1}\kappa_{p,q}^{\vec a, \vec b , 2}\, x^{p-1}y^{q-1}
&=
\sum_{m,n \ge 1} \Bigg[\kappa^{a_1}_{n,m} x^{n-1}y^{m-1}
\sum_{\substack{
   \{p_h\}_{h=1}^n \ge 0 \\
   \{q_h\}_{h=1}^m \ge 0
}}\left(p_1 + 1\right)\left(q_1 +1\right)
\prod_{h=1}^{n} \kappa_{b_1(a_1b_1)^{p_{h}}}x^{p_{h}}
\prod_{h=1}^{m} \kappa_{b_2(a_2b_2)^{q_{h}}}y^{q_{h}}\Bigg]\\[4pt]
&
=\sum_{m,n \ge 1} \Bigg[\kappa^{a_1}_{n,m} x^{n-1}y^{m-1}
\sum_{\substack{
   p_1 \ge 0 \\
}}\left(p_1 + 1\right)\kappa_{b_1(a_1b_1)^{p_{1}}}x^{p_{1}} 
\prod_{h=2}^{n} \sum_{\substack{
   p_{h} \ge 0 \\
  }}\kappa_{b_1(a_1b_1)^{p_{h}}}x^{p_{h}}
  \\[4pt]
&\qquad\qquad\qquad
\cdot \sum_{\substack{
   q_1 \ge 0 \\
}}\left(q_1 + 1\right)\kappa_{b_2(a_2b_2)^{q_{1}}}y^{q_{1}} 
\prod_{h=2}^{n} \sum_{\substack{
  j _{h} \ge 0 \\
  }}\kappa_{b_2(a_2b_2)^{q_{h}}}y^{q_{h}} \Bigg] \\[4pt]
&
=\sum_{m,n \ge 1} \Bigg[\kappa^{a_1}_{n,m} x^{n-1}y^{m-1}
\Bigg(\sum_{\substack{
   p_1 \ge 0 \\
}}\kappa_{b_1(a_1b_1)^{p_{1}}}x^{p_{1}+1}  \Bigg)^{\prime}
\prod_{h=2}^{n} \sum_{\substack{
   p_{h} \ge 0 \\
  }}\kappa_{b_1(a_1b_1)^{p_{h}}}x^{p_{h}}
  \\[4pt]
&\qquad\qquad\qquad 
\cdot \Bigg(\sum_{\substack{
   q_1 \ge 0 \\
}}\kappa_{b_2(a_2b_2)^{q_{1}}}y^{q_{1}+1} \Bigg)^{\prime}
\prod_{h=2}^{n} \sum_{\substack{
  j _{h} \ge 0 \\
  }}\kappa_{b_2(a_2b_2)^{q_{h}}}y^{q_{h}} \Bigg] \\[4pt]
&
=\sum_{m,n \ge 1} \Bigg[\kappa^{a_1}_{n,m} x^{n-1}y^{m-1}
 \Bigg(x\, C_1^{(1)}(x)\Bigg)^{\prime}
(C_1^{(1)}(x))^{n-1} \Bigg(y\, C_1^{(2)}(y)\Bigg)^{\prime}
(C_1^{(2)}(y))^{m-1} \Bigg]\\[4pt]
&
=\Bigg(x\, C_1^{(1)}(x)\Bigg)^{\prime} \Bigg(y\, C_1^{(1)}(y)\Bigg)^{\prime} \sum_{m,n \ge 1} \Bigg[\kappa^{a_1}_{n,m} (x\,C_1^{(1)}(x))^{n-1} 
(y\,C_1^{(2)}(y))^{m-1} \Bigg],
\end{aligned}
\] thereby establishing the desired expression.
\end{proof}

Setting $a = a_1 = a_2$ and $b = b_1 = b_2$ in Lemma~\ref{l2}, $\kappa_{p,q}^{\vec a, \vec b, 2}$ simplifies to  
\begin{equation}
\label{comb}
\begin{aligned}
\kappa^{ab,2}_{p,q}
&=
\sum_{\substack{
\pi=\pi_1\times\pi_2\\
\pi_1\in S_{NC}(p),\ \pi_2\in S_{NC}(q)
}}
\ \sum_{\substack{U_1\in\pi_1\\ U_2\in\pi_2}}
\kappa_{|U_1|,|U_2|}^{a}
\kappa^{ a}_{\pi\setminus\{U_1,U_2\}}
\kappa^{ b}_{Kr(\pi_1) \times Kr(\pi_2)}, 
\end{aligned}
\end{equation}
and in that case we use the notation $R_{ab, 2}$ instead of $R_{\vec a, \vec b, 2}$ for the generating function associated with
\eqref{comb}, which satisfies
\begin{equation}
\label{comb1}
R_{ab,2}(x,y)
=
\big(xC_1(x)\big)'
\big(yC_1(y)\big)'
R_a\big(xC_1(x),yC_1(y)\big).
\end{equation}

The result is very similar to Proposition 6.2 of \cite{CMSS07}, which provides the second contribution to the right-hand side of \eqref{eq3}. The main difference lies in the generating functions appearing as arguments of the two-variable generating function. In the second contribution to \eqref{eq3}, which is the functional reformulation  of the second contribution to \eqref{comb2},  these are $G_a(x)$ and $G_a(y)$. In the present setting, they are replaced by $xC_1(x)$ and $yC_1(y)$, respectively, whose expressions are given in \eqref{eq:C1C2}. The resulting functional equation \eqref{comb1} corresponds to the generalization \eqref{comb}, where non-trivial weights are assigned to the cycles of the Kreweras complements of $\pi_1$ and $\pi_2$.

Combining Lemma~\ref{l1} and \eqref{comb1} we obtain the main theorem of the paper.

\

\section{Applications}\label{si5}

In this section, we shall further assume that $\mathcal{A}$ is an involutive algebra, i.e. there is a conjugate linear map $a \mapsto a^*$ such that $(ab)^* = b^*a^*$. Furthermore, we shall assume that $\varphi(a^*a) \ge 0$ for $a \in \mathcal{A}$.  Such a triad $(\mathcal{A}, \varphi,\varphi_2)$ is frequently called a second order \emph{non-commutative $*$-probability space}.  In the framework of a second order \emph{non-commutative $*$-probability space} we can talk about: selfadjoint random variables, i.e. elements $a \in \mathcal{A}$ with the property that $a = a^*$.

Moreover, throughout these section, the following definitions will be applied.

\begin{defn}[Standard Semicircular  operator]
    A  random selfadjoint variable $s$ in a second order non-commutative $*
$-probability space is called a second-order standard semicircular operator if its first and second order free cumulants satisfy
\[
\label{E}
    \kappa_n^s = \left\{
    \begin{array}{ll}
       1 & \textrm{ if } n=2 \\
        0&  \textrm{ otherwise }
    \end{array}
\right. \; , \qquad \textrm{and} \qquad  \kappa_{p,q}^s =0\;,\ \ \textrm{ for all } p,q\;.
\]
\end{defn}

\begin{defn}[Circular operator]
  Consider $s_{1}$ and $s_
    {2}$ second-order free standard semicircular operators. We call $c =\frac{s_1+is_2}{\sqrt{2}}$
a second-order circular
operator.  
\end{defn}

\subsection{Product with a second-order free Poisson.} \label{as1}

\begin{notation}

Given $a,b$ two second-order non-commutative random variable, we denote by $a\boxplus b$  the second-order non-commutative random variable whose first- and second-order free cumulants are given by
$$
\kappa_p^{a{\boxplus}b}
=
\kappa_p^{a}+\kappa_p^{b}
, \quad \kappa_{p,q}^{a{\boxplus}b}
=
\kappa_{p,q}^{a}+\kappa_{p,q}^{b}\qquad p,q\geq 1.$$

We say that a second-order  non-commutative random variable $b$ is a second-order free Poisson element of rate $t>0$ if
\begin{equation}
\label{eq:free-poisson-rate-t}
\kappa_p^{b}=t,
\qquad
\kappa_{p,q}^b=0,
\qquad p,q\geq 1.
\end{equation}
Given a non-commutative random variable $a$ and $t>0$, let
$$
a_t=ta^{\boxplus 1/t},
$$
where $a^{\boxplus 1/t}$ denotes the non-commutative random variable whose first- and second-order free cumulants are given by
$$
\kappa_p^{a^{\boxplus 1/t}}
=
\frac{\kappa_p^a}{t},
\qquad
\kappa_{p,q}^{a^{\boxplus 1/t}}
=
\frac{\kappa_{p,q}^a}{t}.
$$

We then let $a_{t,1,0}$ denote the non-commutative random variable whose first- and second-order free cumulants satisfy
$$
\kappa_p^{a_{t,1,0}}
=
\kappa_p^{a_t},
\qquad
\kappa_{p,q}^{a_{t,1,0}}
=
0,
$$ so that $\varphi((a_{t,1,0})^p)=\varphi(a_t^p)$. Furthermore, let $\bar{a}_{t,1,1}$ be the non-commutative random variable whose first- and second-order free cumulants are determined by the moments of $a_t$, namely,
$$
\kappa_p^{\bar{a}_{t,1,1}}
=
\varphi(a_t^p),
\qquad
\kappa_{p,q}^{\bar{a}_{t,1,1}}
= \varphi_2(a_t^p,a_t^q).
$$
 Similarly, let $\bar{a}_{t,0,1}$ be the non-commutative random variable whose first- and second-order free cumulants are given by
$$
\kappa_p^{\bar{a}_{t,0,1}}
=0,
\qquad
\kappa_{p,q}^{\bar{a}_{t,0,1}}
=
\varphi_2\bigl((a_{t,1,0})^p,(a_{t,1,0})^q\bigr).
$$

\end{notation}
\begin{rem}
We do not claim the existence of non-commutative random variables with
such prescribed cumulants. The purpose of this notation is merely to
facilitate algebraic manipulations.
\end{rem}
\begin{prop}\label{dap}
Let $b$ be a second-order free Poisson element of rate  $t >0$. Assume, moreover, that the variables $b$ and $a$ are second-order free. Then
\begin{equation*}
{zR_{ab}(z)
=
t\Bigl(\frac{1}{z}G_{a_t}(1/z)-1 \Bigr) \qquad \textrm{ and }\qquad }
R_{ab}(x,y)
=
\frac{1}{x^2y^2}
\Bigl(
tG_{a_t}(1/x,1/y)
+
(1-t)G_{a_{t,1,0}}(1/x,1/y)
\Bigr).
\end{equation*} 
or equivalently for any $p, q\ge 1$:
\begin{equation*}
\kappa_{p}^{ab}
=
t\,\varphi(a_t^p)
 \qquad \text{ and } \qquad
\kappa_{p,q}^{ab}
=
t\,\varphi_2\Bigl(a_t^p,a_t^q\Bigr)
+
(1-t)\,\varphi_2\Bigl((a_{t,1,0})^p,(a_{t,1,0})^q\Bigr).
\end{equation*}
Equivalently, the second order distribution of $ab$ is the same as the distribution of the sum 
$
\bar{a}_{t,1,1}^{\boxplus t}
+
(\bar{a}_{t,0,1})^{\boxplus(1-t)},
$
where $\bar{a}_{t,1,1}^{\boxplus t}$ and $(\bar{a}_{t,0,1})^{\boxplus(1-t)}$ are second order free.

\end{prop}
\begin{proof}
Since $\kappa^{b}_{p,q}=0$ for all $p,q\geq 1$, Theorem~\ref{t1} takes the following form:
$$
R_{ab}(x,y)
= - \frac{\partial^2}{\partial x\,\partial y}
\Bigl[
\log\bigl(
1-xyA(x,y)B(x,y)
\bigr)
\Bigr]
+
\bigl(xC_1(x)\bigr)'\,
\bigl(yC_1(y)\bigr)'\,
R_a\bigl(xC_1(x),yC_1(y)\bigr),
$$
where
$$
A(x,y)
=
\frac{R_a(xC_1(x))-R_a(yC_1(y))}
{xC_1(x)-yC_1(y)}
\qquad\text{and}\qquad
B(x,y)
=
\frac{R_{b}(xC_2(x))-R_{b}(yC_2(y))}
{xC_2(x)-yC_2(y)}.
$$ Since $\kappa_p^{b}=t$ for every $p\geq 1$, we have that
$$
R_{b}(z)=\frac{t}{1-z}.
$$ Consequently,
\begin{equation}
\label{eq:ratio-C2}
\begin{aligned}
B(x,y)
&=
\left(
\frac{t}{1-xC_2(x)}
-
\frac{t}{1-yC_2(y)}
\right)
\frac{1}{xC_2(x)-yC_2(y)}=\frac{t}
{(1-xC_2(x))(1-yC_2(y))}.
\end{aligned}
\end{equation} If $\pi\in S_{NC}(n)$, then, by Definition~\ref{s1d1},
$$
\#(\pi)+\#(\pi^{-1}\gamma_n)=n+1.
$$
Moreover, by \eqref{eq:def-kreweras} and the definition of $Kr^{+1}$,
$$
\#(Kr^{+1}(\pi))
=
\#(Kr(\pi))
=
\#(\pi^{-1}\gamma_n).
$$
Therefore, applying Proposition~\ref{p1}~$(e)$, we obtain, for $n\geq 1$,
$$
\begin{aligned}
\kappa_{b(ab)^n}
&=
\sum_{\pi\in S_{NC}(n)}
\kappa^a_\pi\,
t^{\#(Kr^{+1})(\pi)}=
\sum_{\pi\in S_{NC}(n)}
\kappa^a_\pi\,
t^{n+1-\#(\pi)}=
t^{n+1}
\sum_{\pi\in S_{NC}(n)}
\kappa^{a^{\boxplus 1/t}}_\pi=
t\,\varphi\bigl(a_t^n\bigr).
\end{aligned}
$$ Hence,
\begin{equation}
\label{eq:C1-exemple}
C_1(z)
=
\frac{t}{z}
\sum_{n\geq 0}
\varphi(a_t^n)\,z^{n+1}
=
\frac{t}{z}G_{a_t}(1/z)
=
tM_{a_t}(z).
\end{equation} We also note that
\begin{equation}
\label{eq:Gbox}
tG_{a_t}(1/z)
=
G_{a^{\boxplus 1/t}}(1/tz).
\end{equation} In the same way as above, one has
$$
\kappa_n^{ab}
=
\sum_{\pi\in S_{NC}(n)}
\kappa^a_\pi\,
t^{\#(Kr(\pi))}
=
t\,\varphi(a_t^n).
$$
Therefore,
$$
zR_{ab}(z)
=
\frac{t}{z}G_{a_t}(1/z)-t
=
t\bigl(M_{a_t}(z)-1\bigr),
$$
where we used \eqref{eq2} in the second equality. By Proposition~\ref{p2}~$a)$, $R_{ab}(z)=C_1(z)C_2(z)$, and hence
\begin{equation}
\label{mae}
zC_2(z)
=
1-\frac{1}{M_{a_t}(z)}.
\end{equation}
Thus, \eqref{eq:ratio-C2} becomes
\begin{equation*}
\label{eq:B}
B(x,y)
=
tM_{a_t}(x)M_{a_t}(y)
=
t\frac{G_{a_t}(1/x)G_{a_t}(1/y)}{xy}.
\end{equation*} On the other hand, by \eqref{eq:C1-exemple}, \eqref{eq:Gbox}, and \eqref{eq1}, together with
$$
R_a(z)=tR_{a^{\boxplus 1/t}}(z),
$$
we have
\begin{equation*}
\label{eq:A}
\begin{aligned}
A(x,y)
&=
\frac{
tR_{a^{\boxplus 1/t}}
\bigl(G_{a^{\boxplus 1/t}}(1/tx)\bigr)
-
tR_{a^{\boxplus 1/t}}
\bigl(G_{a^{\boxplus 1/t}}(1/ty)\bigr)
}{
G_{a^{\boxplus 1/t}}(1/tx)
-
G_{a^{\boxplus 1/t}}(1/ty)
}=
\frac{
\left(
\frac{1}{x}
-
\frac{t}{G_{a^{\boxplus 1/t}}(1/tx)}
\right)
-
\left(
\frac{1}{y}
-
\frac{t}{G_{a^{\boxplus 1/t}}(1/ty)}
\right)
}{
G_{a^{\boxplus 1/t}}(1/tx)
-
G_{a^{\boxplus 1/t}}(1/ty)
}.
\end{aligned}
\end{equation*}
We may rewrite this as
$$
A(x,y)
=
\frac{1}{t}
\left(
\frac{1}
{G_{a_t}(1/x)G_{a_t}(1/y)}
+
\frac{\frac{1}{x}-\frac{1}{y}}
{G_{a_t}(1/x)-G_{a_t}(1/y)}
\right).
$$ Combining the expressions above, we obtain
$$
xyA(x,y)B(x,y)
=
1+
G_{a_t}(1/x)G_{a_t}(1/y)
\frac{\frac{1}{x}-\frac{1}{y}}
{G_{a_t}(1/x)-G_{a_t}(1/y)}.
$$ It follows that
$$
\begin{aligned}
R_{ab}(x,y)
={}&
-\frac{\partial^2}{\partial x\,\partial y}
\log\left(
\frac{\frac{1}{x}-\frac{1}{y}}
{G_{a_t}(1/x)-G_{a_t}(1/y)}
\right)+
\frac{\partial}{\partial x}
\bigl(tG_{a_t}(1/x)\bigr)
\frac{\partial}{\partial y}
\bigl(tG_{a_t}(1/y)\bigr)
R_a\Bigl(
tG_{a_t}(1/x),
tG_{a_t}(1/y)
\Bigr).
\end{aligned}
$$ Since
$
\kappa_{n,m}^{a^{\boxplus 1/t}}
=
\frac{\kappa_{n,m}^a}{t},
$
we have
\begin{align*}
R_a\Bigl(
tG_{a_t}(1/x),
tG_{a_t}(1/y)
\Bigr)
=
tR_{a^{\boxplus 1/t}}
\Bigl(
G_{a^{\boxplus 1/t}}(1/tx),
G_{a^{\boxplus 1/t}}(1/ty)
\Bigr).
\end{align*} Moreover,
$$
G_{a^{\boxplus 1/t}}(1/tx,1/ty)
=
t^2G_{a_t}(1/x,1/y).
$$
Hence, combining the preceding relations with \eqref{eq3}, and setting
$
X=\frac{1}{tx}$,
 and 
$Y=\frac{1}{ty},
$
we obtain
\begin{align*}
&\frac{\partial}{\partial x}
\left(tG_{a_t}(1/x)\right)
\frac{\partial}{\partial y}
\left(tG_{a_t}(1/y)\right)
R_a\Bigl(
tG_{a_t}(1/x),
tG_{a_t}(1/y)
\Bigr)
\\
&\qquad=
\frac{1}{tx^2y^2}
\frac{\partial}{\partial X}
\left(G_{a^{\boxplus 1/t}}(X)\right)
\frac{\partial}{\partial Y}
\left(G_{a^{\boxplus 1/t}}(Y)\right)
R_{a^{\boxplus 1/t}}
\Bigl(
G_{a^{\boxplus 1/t}}(X),
G_{a^{\boxplus 1/t}}(Y)
\Bigr)
\\
&\qquad=
\frac{t}{x^2y^2}G_{a_t}(1/x,1/y)
-
t\frac{\partial^2}{\partial x\,\partial y}
\log\left(
\frac{
G_{a_t}(1/x)-G_{a_t}(1/y)
}{
\frac{1}{x}-\frac{1}{y}
}
\right).
\end{align*}
Therefore,
\begin{align*}
R_{ab}(x,y)
={}&
\frac{t}{x^2y^2}G_{a_t}(1/x,1/y)+
(1-t)\frac{\partial^2}{\partial x\,\partial y}
\log\left(
\frac{
G_{a_t}(1/x)-G_{a_t}(1/y)
}{
\frac{1}{x}-\frac{1}{y}
}
\right).
\end{align*} On the other hand,

\begin{equation}
\begin{aligned}
\frac{1}{x^2y^2}G_{a_{t,1,0}}(1/x,1/y)
=
\frac{\partial^2}{\partial x\,\partial y}
\left[
\log\left(
\frac{
G_{a_t}(1/x)-G_{a_t}(1/y)
}{
1/x-1/y
}
\right)
\right].
\end{aligned}
\end{equation}
Consequently,
\begin{equation}
\label{eq:R-transform-rectangular}
R_{ab}(x,y)
=
\frac{1}{x^2y^2}
\Bigl(
tG_{a_t}(1/x,1/y)
+
(1-t)G_{a_{t,1,0}}(1/x,1/y)
\Bigr),
\end{equation}
namely, \begin{equation}
\label{eq:free-cumulants-rectangular}
\kappa_{p,q}^{ab} = t\varphi\bigl(a_t^p, a_t^q\bigr) + (1-t)  \varphi_2\bigl((a_{t,1,0})^p, (a_{t,1,0})^q\bigr). 
\end{equation}
Furthermore since,
$$
R_{\bar{a}_{t,1,1}}(x,y)
=
\frac{1}{x^2y^2}
G_{a_t}(1/x,1/y)
\qquad\text{and}\qquad
R_{\bar{a}_{t,0,1}}(x,y)
=
\frac{1}{x^2y^2}
G_{a_{t,1,0}}(1/x,1/y),
$$ 
it follows that
$$
R_{ab}(z) = t R_{\bar{a}_{t,1,1}}(z) + (1-t)R_{\bar{a}_{t,0,1}}(z)\qquad \textrm{ and } \qquad R_{ab}(x, y) = t R_{\bar{a}_{t,1,1}}(x, y) + (1-t)R_{\bar{a}_{t,0,1}}(x, y).$$
The right-hand side of the above equations are $R$-transforms of the random variable
$$\bar{a}_{t,1,1}^{\boxplus t}
+
(\bar{a}_{t,0,1})^{\boxplus(1-t)},
$$ where the terms in the sum are second order free. This completes the proof.

\end{proof}

\begin{rem}
When $t=1$, this relation simplifies to
$$
zR_{ab}(z)
=
\frac{1}{z}G_{a}(1/z)-1 \qquad \textrm{ and }\qquad 
R_{ab}(x,y)
=
\frac{1}{x^2y^2}\,G_a(1/x,1/y),
$$ or said otherwise, for any $p,q\ge 1$
\begin{equation*}
\kappa_{p}^{ab}
= \kappa_{p}^{cac^*} =
\varphi(a^p)
 \qquad \text{ and } \qquad
\kappa_{p,q}^{ab}
=  \kappa_{p,q}^{cac^*}= 
\varphi_2(a^p,a^q)
,
\end{equation*}
where $c$ is a second-order circular operator, so that $b=cc^*$ is a second-order free Poisson element \eqref{eq:free-poisson-rate-t} of rate $1$. 
 This is consistent with Theorem 8.11 in \cite{AM23}, which states that conjugation by a circular element transforms second-order free cumulants into second-order moments. This observation motivates the discussion in the following section.

In particular, when $a=1$, we obtain
$$
\kappa_{p,q}^{cc^*}=0,
\qquad p,q\geq 1.
$$
\end{rem}

\subsection{Conjugation by a free semicircular element: $sas$.} \label{as2}

\begin{prop} \label{dap1}
    Let $s$ be a second-order standard  semicircular operator. Assume, moreover, that the variables $s$ and $a$ are second-order free. Then the following holds:
\begin{align*} 
R_{as^2}(x,y)
&=
R_{acc^*}(x,y)
-\frac{M_a'(x)M_a'(y)}
{\big(M_a(x)+M_a(y)-M_a(x)M_a(y)\big)^2}, 
\end{align*}
where $c$ is a second-order circular operator.
Consequently,  $$
\kappa_{p,q}^{sas}=\varphi_2(a^p,a^q)+pq \sum_{j=1}^{\min\{p,q\}}\frac{f^a_j(p) f^a_{j}(q)}{j},\qquad p,q\ge1,
$$ where  $$f^a_j(k)=[z^{k}]\bigl(1-M_a(z)\bigr)^{j}=\sum_{i=1}^{j}\binom{j}{i}(-1)^i
\Bigg(\sum_{\substack{
   \{l_h\}_{w=1}^i \ge 0 \\
  l_1+\cdots+l_i=k \\
 }}\ \prod_{w=1}^{i}\varphi\big(a^{l_w}\big)\Bigg).$$
$$$$
\end{prop}
 \begin{proof} It follows from
Theorem~\ref{t1} that
\begin{align*}
R_{as^2}(x,y)
&=
-\frac{\partial^2}{\partial x\,\partial y}
\left[
\log\left(
1-xyA(x,y)B(x,y)
\right)
\right] 
+\big(xC_1(x)\big)'\big(yC_1(y)\big)'
R_a\big(xC_1(x),yC_1(y)\big) \\
&\quad
+\big(xC_2(x)\big)'\big(yC_2(y)\big)'
R_{s^2}\big(xC_2(x),yC_2(y)\big),
\end{align*}
where
\begin{align*}
A(x,y)
&=
\frac{
R_a\big(xC_1(x)\big)-R_a\big(yC_1(y)\big)
}{
xC_1(x)-yC_1(y)
} \quad \text{and }  \quad 
B(x,y)=
\frac{
R_{s^2}\big(xC_2(x)\big)-R_{s^2}\big(yC_2(y)\big)
}{
xC_2(x)-yC_2(y)
}.
\end{align*} It is well-known that (see Lectures 12 and 15 from \cite{NS06}),  in the square case corresponding to \(t=1\), \(s^2\) and \(cc^*\) have the same first-order
distribution; namely, both are free Poisson elements of rate \(1\) (see \eqref{eq:free-poisson-rate-t} at first order).
Consequently, since the first two terms in the above expression
involve only first-order information concerning $s^2$, it follows from Proposition \ref{dap} (for $t=1$) 
that, 
\begin{align*}
R_{as^2}(x,y)
&=R_{acc^*}(x,y) 
 +\big(xC_2(x)\big)'\big(yC_2(y)\big)'
R_{s^2}\big(xC_2(x),yC_2(y)\big),
\end{align*} where $R_{acc^*}(x,y)= \frac{1}{x^2y^2}\,
G_a(1/x,1/y)$. Moreover, according to \eqref{mae},
\begin{equation*}
zC_2(z)
=
1-\frac{1}{M_a(z)}.
\end{equation*} Let us focus on the second term. By Theorem~31 of
\cite{MST09}, we have
\begin{equation*}
R_{s^2}(z,w)
=
\frac{1}{(1-z-w)^2}.
\end{equation*}
Therefore,
\begin{align*}
&\big(xC_2(x)\big)'\big(yC_2(y)\big)'
R_{s^2}\big(xC_2(x),yC_2(y)\big)=
\left(1-\frac{1}{M_a(x)}\right)'
\left(1-\frac{1}{M_a(y)}\right)'
\frac{1}{
\left[
1-\left(1-\frac{1}{M_a(x)}\right)
-\left(1-\frac{1}{M_a(y)}\right)
\right]^2
} \\
&\qquad=
\frac{M_a'(x)M_a'(y)}
{\big(M_a(x)+M_a(y)-M_a(x)M_a(y)\big)^2} =
-\frac{\partial^2}{\partial x\,\partial y}
\left[
\log\left(
M_a(x)+M_a(y)-M_a(x)M_a(y)
\right)
\right],
\end{align*} which proves the first statement. Now, note that
$$
-\log\Big(M_a(x)+M_a(y)-M_a(x)M_a(y)\Big)
= -\log\Big(1-F_a(x)F_a(y)\Big)
= \sum_{j\geq 1}\frac{F_a(x)^{j}F_a(y)^{j}}{j},
$$
where $F_a(z):=1-M_a(z)$. 
Letting $F_a(x,y):=-\frac{\partial^2}{\partial x\,\partial y}
\log\Big(1-F_a(x)F_a(y)\Big)$, we therefore have that:
$$
[x^{p-1}y^{q-1}]F_a(x,y) = - pq [x^{p}y^{q}] \log\Big(1-F_a(x)F_a(y)\Big) = pq \sum_{j\geq 1}\frac{f^a_j(p)\,f^a_j(q)}{j}\;,
$$ where $f^a_j(k):=[z^{k}]F_a(z)^{j}$. Note that since $F_a(0)=1-M_a(0)=0$, we have that $f^a_j(k)=0$ whenever $j>k$.

Finally for $1 \leq j \leq k$ we expand the binomial:
$$
\big(1-M_a(z)\big)^{j}=1+\sum_{i=1}^{j}\binom{j}{i}(-1)^i M_a(z)^i,
$$  
so that
$$
f^a_j(k)=\sum_{i=1}^{j}\binom{j}{i}(-1)^i\,[z^{k}]M_a(z)^i
=\sum_{i=1}^{j}\binom{j}{i}(-1)^i
\Bigg(\sum_{\substack{
   \{l_w\}_{w=1}^i \ge 0 \\
  l_1+\cdots+l_i=k \\
 }}\ \prod_{w=1}^{i}\varphi\big(a^{l_w}\big)\Bigg).
$$
 which completes the proof.

 \end{proof}

\begin{rem}
Propositions \ref{dap} and \ref{dap1} show that, although $sas$ and $cac^*$ have the same first-order distribution, as established in Proposition 12.18 of \cite{NS06} and Example 8.10 of \cite{AM23}, respectively, their second-order distributions are different. Indeed, the second-order cumulants of $sas$ contain an additional contribution. 

An important special case of the preceding propositions is obtained by setting $a=1$, in which case $sas=s^2$ and $cac^*=cc^*$. In this setting, Proposition \ref{dap1} is consistent with Proposition 30 of \cite{MST09}. Moreover, Examples 8.4 and 8.6 of \cite{ADS25} show that $s^2$ and $cc^*$ have the same third-order distribution. Hence, although $s^2$ and $cc^*$ agree at the first and third orders, they exhibit distinct behavior at the second order. This provides another illustration of the subtle distinctions that can arise among higher-order distributions in non-commutative probability spaces.

\end{rem}

\appendix

\section{An Application to Finite Free Probability}\label{si6}

One big success of free probability is its application in describing the asymptotic spectral behavior of large random matrices: independent, unitarily invariant matrices become asymptotically free as their size tends to infinity, and operations such as the free additive and multiplicative convolutions describe the limiting eigenvalue distributions of sums and products.

Finite free probability is a more recent framework that transfers these ideas to the realm of polynomials. Instead of working with random matrices, it works with polynomials of fixed degree $d$, identified with the characteristic polynomials of $d \times d$ matrices. Its main objects of interest are the finite free convolutions, which are bilinear operations on polynomials that preserve positivity properties, i.e. real rootness. 

The original motivation of the foundational paper of Marcus Spielman and Srivastava~\cite{MSS2022} was to obtain bounds for graph constructions used in the proof of existence of Ramanujan graphs.
However, the connection of this work with free probability, made precise later,  has since
proved fruitful in several directions. It has found applications to the theory of orthogonal polynomials~\cite{martinez2026zeros, jalowy2026zeros}, to Appell
sequences~\cite{campbell2026free, campbell2026p}, to entire functions~\cite{campbell2024universality},
and to the behaviour of roots under the heat flow~\cite{hall2023heat}. It is also
related to Harish-Chandra--Itzykson--Zuber (HCIZ) integrals \cite{mergny2021rank, benaych2022matrix}, and has in turn given
new insights into free probability itself~\cite{arizmendi2026free, arizmendi2026free1}.

Here we are interested in finite free multiplicative convolution. Let $\mathbb{R}^+_d[x]$ denote the set of all polynomials in $x$ of degree $d$ with non-negative roots.
Given any monic polynomials $p_d,q_d\in \mathbb{R}^+_d[x] $ the \emph{finite-free multiplicative convolution} of $p_d$ and $q_d$ is
defined as the polynomial $p_d\boxtimes_d q_d \in \mathbb{R}^+_d[x]$ given by 
\begin{equation}\label{eq:def.multiplicative}
[p_d \boxtimes_d q_d](x) =
			\sum_{k=0}^d x^{d-k}(-1)^k \binom{d}{k} \Big[ \coef{k}{p_d} \coef{k}{q_d} \Big],
\end{equation}
where, for a polynomial of degree $d$, $\coef{k}{r_d}$ denotes its $k$-th normalized coefficients, determined by the equation $r_d (x) =\sum_{k=0}^d x^{d-k}(-1)^k \binom{d}{k}  \coef{k}{r_d}$.

            This operation was studied independently by Szeg\"{o}~\cite{szego1922bemerkungen}
and Walsh~\cite{walsh1922location}.  It was rediscovered by Marcus,
Spielman, and Srivastava~\cite{MSS2022} as the expected characteristic polynomial
of a product of randomly rotated non-negative matrices.  

A fundamental result proved in Arizmendi, Garza-Vargas and Perales
\cite{APV25} is that finite free multiplicative convolution approximates free convolution as the degree tends to infinity. More precisely, for a polynomial $p_d$ of degree $d$ with $d$ real roots  counted with multiplicities,  say 
$
\lambda_1(p_d),\ldots,\lambda_d(p_d)$, we denote by $\mu_{p_d}$,  its empirical root distribution, given by 
\[
\mu_{p_d}
=
\frac1d\sum_{j=1}^d\delta_{\lambda_j(p_d)}.
\]
Suppose that two sequences of polynomials $\mfp=(p_d)_{d\geq1}$ and $\mfq=(q_d)_{d\geq1}$ have limiting root distributions in the sense that 
$\mu_{p_d}\to\mu,$ and $\mu_{q_d}\to\nu$
in distribution, as $d\to \infty$. Then 
$$\mu_{p_d\boxtimes_d q_d} \to\mu \boxtimes\nu.$$

We are interested in the $1/d$ correction of the finite free multiplicative convolution, $\boxtimes_d$  as \( d \) goes to infinity, following the work of \cite{APV25}. This can be formulated, when $\mu$ and $\nu$ have compact support, as the existence, for each $k\in \mathbb{N}$, of the following limit 
$$m_k'(\mfp\boxtimes\mfq) =\lim_{d\to \infty} \left(\int t^k d\mu_{p_d\boxtimes_d q_d}-\int t^k d\mu \boxtimes\nu\right).$$
In this case, we call the sequence $\{m_k'(\mfp\boxtimes\mfq)\}_{k\geq1}$, the \emph{infinitesimal moments} of the sequence $\mfp\boxtimes\mfq:=(p_d\boxtimes_d q_d)_{d\geq1}$.
We refer the reader to \cite{AGVP23} and \cite{APV25}  for further details on this problem and the related results.  

More generally when a sequence of polynomials $\mfp= (p_d)_{d\geq1}$ has infinitesimal moments of all orders, we identify the sequence $(m_n'(\mfp))_{n\geq1}$ with the unique linear functional $\mu':\rr[x] \to \rr$ such that $\mu'(x^n) = m_n'(\mfp)$ for every $n \geq 0$. 
Note that we always have $\mu'(1)=0$. 
We refer to the pair $(\mu,\mu')$ as the \emph{infinitesimal distribution} of $\mfp = (p_d)_{d\geq1}$. 
We also use the notation $m_n(\mu')=m_n'(\mfp)$. The associated infinitesimal Cauchy transform is defined as: 
\[
G_{\mu'}(z) =\sum_{n\ge 1} m_n(\mu') z^{-(n+1)}.
\]

The notion of infinitesimal distribution $(\mu,\mu')$, for a non-commutative variable, appeared first in the work of Belinschi and Shlyakhtenko \cite{belinschi2012free}, as a generalization of the Type B distributions of Biane, Goodman, and Nica \cite{biane2003non} and in this work they only require $\mu$ and $\mu'$ to be linear functionals on $\mathbb{C}[x]$ with $\mu(1)=1$ and $\mu'(1)=0$.

Throughout the appendix, we consider two sequences of polynomials $\mfp = (p_d)_{d\geq1}$ and $\mfq = (q_d)_{d\geq1}$ with infinitesimal distributions $(\mu,\mu')$ and $(\nu,\nu')$, respectively. Moreover, $\kappa^{(d)}_n$ denotes the finite free cumulants and $\widehat{r}_n$ the infinitesimal part
of the finite free cumulants. 

By Lemma~3.1 from \cite{AGVP23}, the hypothesis of having infinitesimal moments is equivalent, to the existence of measures $\mu$ and $\nu$ and sequences $(\widehat r_n(\mfp))_{n\ge1}$ and $(\widehat r_n(\mfq))_{n\ge1}$ such that for all $n\in\mathbb{N}$:
\begin{equation} \label{inf cumulants}
\kappa^{(d)}_n(p_d)= \kappa_n(\mu) + \frac{1}{d}\,\widehat r_n(\mfp) + o(1/d)
\qquad\text{and}\qquad
\kappa^{(d)}_n(q_d)= \kappa_n(\nu) + \frac{1}{d}\,\widehat r_n(\mfq) + o(1/d).
\end{equation}

Under the above conditions, according to Theorem 1.1 from  \cite{AGVP23}, the sequence of finite free multiplicative convolutions, $\mfp \boxtimes \mfq := (p_d \boxtimes q_d)_{d=1}^\infty$ 
has as limit distribution $\mu \boxtimes \nu$ and the infinitesimal distribution satisfies that for all $n\in\mathbb{N}$, the finite free cumulant $\kappa^{(d)}_n(p_d \boxtimes q_d)$ has the asymptotic expansion 
\begin{equation*}
\kappa^{(d)}_n(p_d \boxtimes q_d)
=
\kappa_n(\mu \boxtimes \nu)
+ \frac{1}{d}\,\widehat r_n(\mfp\boxtimes\mfq)
+ o(1/d),
\end{equation*} 
where
\begin{align}\label{infinitesimalproduct}
\widehat r_n(\mfp\boxtimes\mfq)
&=
\sum_{\pi\in S_{NC}(n)}
\Bigg(
    \sum_{V\in\pi}
        \widehat r_{|V|}(\mfp)\, \kappa_{\pi\setminus V}(\mu)\, \kappa_{Kr(\pi)}(\nu)
    + \sum_{W \in Kr(\pi)}
        \widehat r_{|W|}(\mfq)\, \kappa_{\pi}(\mu)\, \kappa_{Kr(\pi)\setminus W}(\nu)
\Bigg)
\\[1ex]
&\quad
- \frac{n}{2}
\sum_{\substack{t+s=n \\ \pi\in S_{NC}(t,s)}}
\frac{ \kappa_\pi(\mu) \kappa_{Kr(\pi)}(\nu)}{ts}.
\nonumber
\end{align}
For any sequence of polynomials  $\mfp = (p_d)_{d\geq1}$, we let 
\[
\widehat{R}_{\mfp}(z)
=
\sum_{n\geq 1}\widehat{r}_n(\mfp)z^{n-1}. 
\]
Let us state the main result of this appendix, which provides a formula for the generating series $\widehat{R}_{\mfp\boxtimes \mfq}(z)$.

\begin{thm}\label{ffp3}
   Let $\mu$ and $\nu$ be probability measures with finite moments
of all orders. The generating series of the infinitesimal part of the finite cumulants series
of $\mfp\boxtimes\mfq$ is given by
\begin{align*}
\widehat{R}_{\mfp\boxtimes\mfq}(z)
={}&
\big(xC_1^{\mu,\nu}(z)\big)'
\widehat{R}_{\mfp}\big(zC_1^{\mu,\nu}(z)\big)
+
\big(zC_2^{\mu,\nu}(z)\big)'
\widehat{R}_{\mfq}\big(zC_2^{\mu,\nu}(z)\big) 
\\
&\hspace{2cm} +
\frac{1}{2}\frac{d}{dz}
\log\Bigl(
1-z^2
R_\mu'\big(zC_1^{\mu,\nu}(z)\big)
R_\nu'\big(zC_2^{\mu,\nu}(z)\big)
\Bigr),
\end{align*}
where $C_1^{\mu,\nu}(z)$ and  $C_2^{\mu,\nu}(z)$ are determined by the equations
$$C_2^{\mu,\nu}(z)=R_{\mu}(z\,C_1^{\mu,\nu}(z))\quad\text{ and} \quad C_1^{\mu,\nu}(z)=R_{\nu}(z\,C_2^{\mu,\nu}(z)).$$
\end{thm}

There are two kinds of contributions, corresponding to the terms on the first and second line of the statement. We will prove these contributions in Lemmas \ref{ffp1} and \ref{ffp2}, respectively. The former is a consequence of Lemma~\ref{l2}: it provides an explicit formula for the first term in the sum of the right-hand side of \eqref{infinitesimalproduct}.  Let us remark that in principle this sum can be handled from the original theory of infinitesimal free probability, for instance by using triangular matrix realization of infinitesimal freeness, see for example Section 1 in \cite{belinschi2012free}. However, we here give a proof from our main formulas to show the connection between the two theories.

\begin{lem}
    \label{ffp1}
   Let $\mu$ and $\nu$ be probability measure with all moments, and $\mfp = (p_d)_{d\geq1}$ as in \eqref{inf cumulants} and define for each $n\in\mathbb{N}$,  
\[
\widehat{r}_n^{(1)}(\mfp\boxtimes\mfq)
= \sum_{\pi\in S_{NC}(n)}
    \sum_{V\in\pi}
        \widehat r_{|V|}(\mfp)\, \kappa_{\pi\setminus V}(\mu)\, \kappa_{Kr(\pi)}(\nu)
\]
then its generating series is given by
\[
\sum_{n \geq 1 }\widehat{r}_n^{(1)}(\mfp\boxtimes\mfq)\, z^{\,n-1}
=  \Big(zC_1^{\mu,\nu}(z)\Big)'
\hat{R}_{\mfp}\Big(zC_1^{\mu,\nu}(z)\Big).
\]
where $C^{\mu,\nu}_1$ is determined by the equation
$$C^{\mu,\nu}_1(z)=R_{\nu}\Big(z\,R_{\mu}\big(z\,C^{\mu,\nu}_1(z)\big)\Big).$$ 
\end{lem}
\begin{proof}
Setting $a_1$ and $b_1$ in Lemma \ref{l2} to be first-order free
random variables satisfying
\[
\kappa_p^{a_1} = \kappa_p^\mu,
\qquad
\kappa_p^{b_1} = \kappa_p^\nu,
\]
and
\[
\kappa_{p,q}^{a_1} =
\begin{cases}
\widehat{r}_p(\mfp), & \text{if } q=1,\\
0, & \text{if } q\geq 2.
\end{cases}
\]Similarly, let $a_2$ and $b_2$ be first-order free random variables
such that
\[
\kappa_p^{a_2} = 0,
\qquad
\kappa_p^{b_2} =
\begin{cases}
1, & \text{if } p=1,\\
0, & \text{if } p\geq 2.
\end{cases}
\] It follows that
\[
\kappa_{p,q}^{\vec a, \vec b , 2} =
\begin{cases}
\widehat{r}_p^{(1)}(\mfp\boxtimes\mfq), & \text{if } q=1,\\
0, & \text{if } q\geq 2.
\end{cases}
\]
Therefore, by Lemma \ref{l2},
\begin{align*}
\sum_{n\geq 1} \widehat{r}_n^{(1)}(\mfp\boxtimes\mfq)\,z^{n-1}
= R_{\vec a, \vec b, 2}(z,y)= \big(zC_1^{(1)}(z)\big)'
   \big(yC_1^{(2)}(y)\big)'
   R_a\big(zC_1^{(1)}(z),yC_1^{(2)}(y)\big).
\end{align*} Observe that
\[
R_a(z,y)=\widehat{R}_{\mfp}(z).
\]
Hence, for every $y$,
\begin{align*}
\sum_{n\geq 1} \widehat{r}_n^{(1)}(\mfp\boxtimes\mfq)\,z^{n-1}
= \big(zC_1^{(1)}(z)\big)'
   \Big(C_1^{(2)}(y)+y(C_1^{(2)}(y))'\Big)
   \widehat{R}_{\mfp}(z).
\end{align*} Finally, setting $y=0$ yields the desired result.
\end{proof}

We now turn to the problem of finding an explicit formula for the sum on the second line of
\eqref{infinitesimalproduct}, which involves a sum over annular
non-crossing partitions in $S_{NC}(t,s)$. This question was left open by
Arizmendi, Perales, and Vázquez-Becerra in
\cite[Section 5.4]{APV25}. To address it, we will use an expression of a two variable series obtained from that of our main theorem by integration, followed by
a careful analysis of the resulting constants. This is the content of the next lemma. 

For $u$ a second-order non-commutative random variable, we let:
$$
T_{u}(z,y)=\sum_{i,j\geq 1} \kappa_{i,j}^u \frac{z^{i}y^{j}}{ij}.
$$

\begin{lem} \label{c1} Let $a$ and $b$ be two second order free non-commutative random variables. Then the following holds:
\begin{align*}
 T_{ab}(z,y)
&= -
 \log\left(
1
-
z\,y \bigg(\frac{R_a(z\,C_1(z))-R_a(y\,C_1(y))}
         {z\,C_1(z)-y\,C_1(y)}\bigg)\,\bigg(\frac{R_b(z\,C_2(z))-R_b(y\,C_2(y))}
     {z\,C_2(z)-y\,C_2(y)}\bigg)
\right) \nonumber\\[4pt]
&\quad + \,
T_a\Big(z\,C_1(z),\,y\,C_1(y)\Big)  + 
T_b\Big(z\,C_2(z),\,y\,C_2(y)\Big). 
\end{align*}
\end{lem} 

\begin{proof}
By definition, we have that

$$\frac{\partial^2}{\partial z\,\partial y} T_u(z,y)=R_u(z,y).$$ Hence, as a direct consequence of Theorem \ref{t1}:

\begin{align*}
\frac{\partial^2}{\partial z\,\partial y} T_{ab}(z,y)
&= -\frac{\partial^2}{\partial z\,\partial y}
\left[ \log\left(
1
-
z\,y \bigg(\frac{R_a(z\,C_1(z))-R_a(y\,C_1(y))}
         {z\,C_1(z)-y\,C_1(y)}\bigg)\,\bigg(\frac{R_b(z\,C_2(z))-R_b(y\,C_2(y))}
     {z\,C_2(z)-y\,C_2(y)}\bigg)
\right) \right]\nonumber\\[4pt]
&\quad + \,
\frac{\partial^2}{\partial z\,\partial y} 
T_a\big(z\,C_1(z),\,y\,C_1(y)\big) \nonumber +
\frac{\partial^2}{\partial z\,\partial y}
T_b\big(z\,C_2(z),\,y\,C_2(y)\big). 
\end{align*}
Thus, 
\begin{align}\label{cr}
 T_{ab}(z,y)
&= -
 \log\left(
1
-
z\,y \bigg(\frac{R_a(z\,C_1(z))-R_a(y\,C_1(y))}
         {z\,C_1(z)-y\,C_1(y)}\bigg)\,\bigg(\frac{R_b(z\,C_2(z))-R_b(y\,C_2(y))}
     {z\,C_2(z)-y\,C_2(y)}\bigg)
\right) \nonumber\\[4pt]
&\quad + \,
T_a\big(z\,C_1(z),\,y\,C_1(y)\big)+ 
T_b\big(z\,C_2(z),\,y\,C_2(y)\big) + 
g_1(z)+g_2(y),
\end{align} Differentiating the above expression with respect to $z$, we get:
\begin{align*}
g_1'(z)  
&= \frac{\partial}{\partial z} \bigg[ T_{ab}(z,y) 
+ \log\left(
1
-
z\,y \bigg(\frac{R_a(z\,C_1(z))-R_a(y\,C_1(y))}
         {z\,C_1(z)-y\,C_1(y)}\bigg)\,\bigg(\frac{R_b(z\,C_2(z))-R_b(y\,C_2(y))}
     {z\,C_2(z)-y\,C_2(y)}\bigg)
\right) \nonumber\\[4pt]
&\quad -\,
T_a\big(z\,C_1(z),\,y\,C_1(y)\big) - \,
T_b\big(z\,C_2(z),\,y\,C_2(y)\big) \bigg].
\end{align*}
Hence, since the above equation is valid for any $y$, letting $y$ approach $0$, it follows that 

\[
g_1'(z)= 0 .
\] Similarly, we obtain $g_2'(y)=0$.  Therefore, $g_1(z)+g_2(y)=c,$ where $c$ is a constant. Finally, letting $y=0$ in \eqref{cr}, we obtain $c=0$, which completes the proof.
\end{proof}

Applying now Lemma~\ref{c1} and analyzing the behavior of the resulting expression as
$y$ approaches $z$ yields the desired formula. This is the 
content of the following lemma.

\begin{lem}\label{ffp2}
   Let $\mu$ and $\nu$ be probability measures with all moments, and define for each $n\in\mathbb{N}$,  
\[
\widehat{r}_n^{(2)}(\mfp\boxtimes\mfq)
= -\frac{n}{2}
\sum_{\substack{t+s=n \\ \pi \in S_{NC}(t,s)}}
\frac{ \kappa_{\pi}(\mu)\, \kappa_{K r(\pi)}(\nu)}{t s}.
\]
Then the generating series for this sequence of numbers is given by
\[
\sum_{n \geq 2}\widehat{r}_n^{(2)}(\mfp\boxtimes\mfq)\, z^{\,n-1}
= \frac{1}{2}\frac{d}{ dz}
\log \bigg(1 -
z^2
    \cdot R_{\mu}'(z C_{1}(z)) \cdot
R'_{\nu}(z C_{2}(z))
\bigg),
\]
where for a measure $\rho$,  $ R_\rho(z)=(C_{\rho}(z)-1)/z$.
\end{lem}

\begin{proof}
    By Lemma \ref{l1} and by the same reasoning as in Lemma \ref{c1},

\begin{equation*} 
\begin{aligned}
-\sum_{t,s\ge 1}
\Bigg(\sum_{\substack{\pi\in S_{NC}(t,s)}}
\kappa_{\pi}(\mu)\, \kappa_{K r(\pi)}(\nu)\Bigg)\, \frac{z^{t} y^{s}
}{ts}
= \log\left(
1
-
z\,y \bigg(\frac{R_{\mu}(z\,C_1(z))-R_{\mu}(y\,C_1(y))}
         {z\,C_1(z)-y\,C_1(y)}\bigg)\,\bigg(\frac{R_{\nu}(z\,C_2(z))-R_{\nu}(y\,C_2(y))}
     {x\,C_2(x)-y\,C_2(y)}\bigg)
\right).
\end{aligned}
\end{equation*}
At this point, note that when $y$ approaches $z$, the left-hand side of the above expression, we get that
\begin{equation}
\label{eq11}
\sum_{n\geq 2} \bigg(-\sum_{\substack{t+s=n \\ \pi \in S_{NC}(t,s)}}
\frac{ \kappa_{\pi}(\mu)\, \kappa_{K r(\pi)}(\nu)}{t s} \bigg) z^n.
\end{equation}
Now letting $y$ approach
$z$ in the right-hand side of the same expression becomes 

\begin{equation}\label{me2}
\begin{aligned}
\log\left( \lim_{y \to z}\left(
1 -
z\,y \bigg(\frac{R_{\mu}(z\,C_1(z))-R_{\mu}(y\,C_1(y))}
         {z\,C_1(z)-y\,C_1(y)}\bigg)\,\bigg(\frac{R_{\nu}(z\,C_2(z))-R_{\nu}(y\,C_2(y))}
     {z\,C_2(z)-y\,C_2(y)}\bigg)
\right)\right) \\
= \log \bigg(1 -
z^2 \cdot R_{\mu}'(z C_{1}(z)) \cdot
R'_{\nu}(z C_{2}(z))
\bigg).
\end{aligned}
\end{equation}
 Finally, after substituting and differentiating with respect to $z$ in \eqref{eq11} and \eqref{me2}, and then dividing by 2, the proof is complete.
\end{proof}

The following result is Theorem 5.5 from \cite{AGVP23}. We will derive it as a particular case of the above proposition.
\begin{cor}\label{dgp}
Suppose that  for any integer $n\ge 1$,
\[
m_n'(\mu)
= -\frac{n}{2}\ 
\sum_{\substack{t,s\ge 1,\ \ t+s=n \\ \pi \in S_{NC}(t,s)}}
\frac{ \kappa_{\pi}(\mu)}{t s},
\]
then
\[
\sum_{n \geq 2}m_n'(\mu)\, z^{\,-(n+1)}
 = \frac{G''_{\mu}(z)}{2\,G'_{\mu}(z)} - \frac{G'_{\mu}(z)}{G_{\mu}(z)} .\]
Moreover, if $\mfp$ has limiting distribution $\mu$ and satisfies that $\widehat{R}_{\mfp}(z)=0,$ then $m_n(\mu')=m_n'(\mu)$ and then the quantity above is $G_{\mu'}(z) $.
\end{cor}
  
\begin{proof}
Consider $\nu$ to be a free Poisson element with rate $1$, i.e. $\kappa_n(\nu)=1$, for all $n$. In this case 
\[
\widehat{r}_n^{(2)}(\mfp\boxtimes\mfq)
= -\frac{n}{2}
\sum_{\substack{t+s=n \\ \pi \in S_{NC}(t,s)}}
\frac{ \kappa_{\pi}(\mu)}{t s}=m_n'(\mu).
\]
By Lemma \ref{ffp2} we have that

$$\sum_{n \geq 2}\widehat{r}_n^{2}(\mfp\boxtimes\mfq)\, z^{\,-(n-1)}=-\frac{1}{2} \cdot z^2 \cdot\frac{d}{ dz}
\log \bigg(1 -
z^{-2}
    \cdot R_{\mu}'(z^{-1} C_{1}(z^{-1})) \cdot
R'_{\nu}(z^{-1} C_{2}(z^{-1}))
\bigg),
$$ and then 
$$\sum_{n \geq 2}m_n'(\mu)\, z^{\,-(n+1)}
 = -\frac{1}{2} \cdot\frac{d}{dz}
\log \bigg(1 -
z^{-2}
    \cdot R_{\mu}'(z^{-1} C_{1}(z^{-1})) \cdot
R'_{\nu}(z^{-1} C_{2}(z^{-1}))
\bigg)$$ 
 
 Now, similar to  Proposition  \ref{dap},
$$R_{\nu}(z)= \frac{1}{1-z}\qquad \textrm{and}\qquad R_{\nu}'(z)= \bigl(R_{\nu}(z)\bigr)^2,$$ so that applying the statement $(e)$ of Prop.~\ref{p2}, followed by \eqref{eq:C1-exemple} (for $t=1$): $$R'_{\nu}(z^{-1} C_{2}(z^{-1}))=C_1(z^{-1})^2=(z\, G_{\mu}(z))^2.$$ 
On the other hand,  from the last equality and differentiating $(\ref{eq1})$,
$$R_{\mu}'(z^{-1} C_{1}(z^{-1})) =R'_{\mu}(G_{\mu}(z))=\frac{1}{G'_{\mu}(z)}+\frac{1}{G^2_{\mu}(z)}.$$ Therefore, after simplification, we obtain
$$1 -z^{-2}
    \cdot R_{\mu}'(z^{-1} C_{1}(z^{-1})) \cdot
R'_{\nu}(z^{-1} C_{2}(z^{-1}))=-\frac{G^2_{\mu}(z)}{{G^{\prime}_{\mu}(z)}}.$$  The proof follows by taking the logarithmic derivative and multiplying by $-1/2$. 
\end{proof}

Notice in the proof of the previous corollary that we calculate $\sum_{n \geq 2}\widehat{r}_n^{(2)}(\mfp\boxtimes\mfq)\, z^{\,n-1}
$ when $\nu$ is a Poisson distribution.  
We revisit this calculation on Example \ref{freepoisson} in relation with Laguerre polynomials.

\subsection*{Examples on finite free multiplicative convolution}

In this subsection we present some applications of Theorem~\ref{ffp3}. We first recall a lemma that explains the relation between the infinitesimal Cauchy transforms and $\widehat{R}_{\mfp}$.

\begin{lem}[{\cite[Lemma~4.3]{APV25}}] \label{Lemma H}
Let $\mfp=(p_d)_{d\geq1}$ be a sequence of polynomials with infinitesimal distribution $(\mu,\mu')$. Then
\begin{equation}\label{qis0}
G_{\mu'}(z) = -\widehat{R}_{\mfp}\big(G_{\mu}(z)\big)\, G'_{\mu}(z) - \gdos_{\mu}(z),
\end{equation}
where
\[
\gdos_{\mu}(z):=\frac{G'_{\mu}(z)}{G_{\mu}(z)}-\frac{G''_{\mu}(z)}{2\,G'_{\mu}(z)}.
\]
\end{lem}

In the next two examples, we consider an important case that was left open in \cite{APV25}. Suppose that the finite free cumulants of $\mfp$ and $\mfq$ coincide with the free cumulants of their limiting distributions $\mu$ and $\nu$, respectively. Then $\hat r_n(\mfp)=\hat r_n(\mfq)=0$ for all $n\geq1$, that is,
\[
\widehat{R}_{\mfp}=\widehat{R}_{\mfq}=0,
\]
and \eqref{qis0} reduces to
\[
G_{\mu'}=-\gdos_\mu \qquad\text{and}\qquad G_{\nu'}=-\gdos_\nu .
\]
Although restrictive, the condition of vanishing infinitesimal finite free cumulants is satisfied by various classical sequences of polynomials, such as the Hermite and Laguerre polynomials; see \cite[Section~5.1]{AGVP23}.

Let $(\rho,\rho')$ denote the infinitesimal distribution of the sequence $\mfp\boxtimes\mfq=(p_d\boxtimes_d q_d)_{d\geq1}$. Here  $\rho=\mu\boxtimes\nu$ and the above assumptions, a direct application of Theorem~\ref{ffp3} and Lemma~\ref{Lemma H} yields
\begin{equation}\label{Grho}
G_{\rho'} = -\big(\widehat{R}_{\mfp\boxtimes\mfq}\circ G_{\rho}\big)\, G'_{\rho}-\gdos_\rho ,
\end{equation}
where
\begin{equation}\label{no inf cum}
\widehat{R}_{\mfp\boxtimes\mfq}(z)=\frac{1}{2}\,\frac{d}{dz}
\log\Big(1-z^2\, R_{\mu}'\big(z\, C_{1}^{\mu,\nu}(z)\big)\, R'_{\nu}\big(z\, C_{2}^{\mu,\nu}(z)\big)\Big).
\end{equation}

\begin{exam}[Multiplication with Laguerre] \label{freepoisson}
Let $\mfq=(q_d)_{d\geq1}$ be sequence of normalized Laguerre polynomials, as in \cite[Example 2.8]{AGVP23}. They are defined by its finite free cumulants $\kappa^{(d)}_n(q_d)=1$, for all $n$, and satisfy $
\widehat{R}_{\mfq}=0$. The limiting distribution of $\mfq$, denoted by $\nu$, is a free Poisson distribution with rate $1$. Assume further that the sequence $\mfp=(p_d)_{d\geq1}$ satisfies 
$\widehat{R}_{\mfp}=0$.
Then as in the proof of Corollary~\ref{dgp}, we observe that the equation $R_{\nu}(z)= \frac{1}{1-z}$, implies that
$$R'_{\nu}(z C_{2}(z))=(1/z\, G_{\mu}(1/z))^2\quad \text{and} \quad R_{\mu}'(z C_{1}(z))=\frac{1}{G'_{\mu}(1/z)}+\frac{1}{G^2_{\mu}(1/z)}.$$ Therefore, we obtain
$$\log\left(1 -z^{2}
    \cdot R_{\mu}'(z C_{1}(z)) \cdot
R'_{\nu}(z C_{2}(z))\right)=\log\left(-\frac{G^2_{\mu}(1/z)}{{G^{\prime}_{\mu}(1/z)}}\right)=2\log G_{\mu}(1/z)-\log(-G_\mu(1/z))$$  After differentiating and dividing by 2, we obtain
\begin{equation}\label{Rhatpoisson}
\widehat{R}_{\mfp\boxtimes\mfq}(z) = -\frac{1}{z^2}\, \gdos_\mu \left(\frac{1}{z}\right).
\end{equation}
\end{exam}

\begin{rem}
Example~\ref{freepoisson} corresponds to \cite[Example~5.14]{APV25}; however, the final expression obtained there was incorrectly identified as $H_\mu(z)$ instead of $-1/z^2H_\mu(1/z^2)$.
\end{rem}

Let us now treat a concrete case, when $\mfp$ is the sequence of Hermite polynomials.

\begin{exam}[Multiplication of Hermite and Laguerre]
We now consider the finite free multiplication of Hermite polynomials and Laguerre. This fits the previous example, which means that $\nu$ is a Poisson distribution with rate $1$,  $\mu$ is  the standard semicircle distribution, see e.g.\cite[Example 2.7]{AGVP23}, and $\widehat{R}_{\mfp}=\widehat{R}_{\mfq}=0$.  We will show that \begin{equation}
\label{eq:sec-ex}
G_{\rho'}(z)=-\,\frac{2\,G_\rho^3(z)\,\Big(2-\sqrt{1-4G_\rho^2(z)}\Big)}
{\Big(1-2\sqrt{1-4G_\rho^2(z)}\Big)^2}.
\end{equation}
Since $G_\mu=\frac{1}{2}(z-\sqrt{z^2-4})$ we see that
\[
\gdos_\mu(z)=\frac{z-\sqrt{z^2-4}}{2\,(z^2-4)} .
\]
Thus, by \eqref{Rhatpoisson},
\begin{equation}\label{Rhatsemicircle}
\widehat{R}_{\mfp\boxtimes\mfq}(z)
=-\frac{1}{z^2}\,\gdos_\mu\left(\frac1z\right)
=-\frac{1-\sqrt{1-4z^2}}{2z\,(1-4z^2)} .
\end{equation} We will use the following auxiliary function to simplify $G_{\rho'}$,
\[
g(z):=\sqrt{1-4\,G_{\rho}^2(z)}.
\] Substituting \eqref{Rhatsemicircle} into \eqref{Grho}, we obtain
\begin{equation}\label{examff}
G_{\rho'}(z)
= \frac{1-g(z)}{g^2(z)}\cdot\frac{G'_{\rho}(z)}{2G_{\rho}(z)}
 - \frac{G'_{\rho}(z)}{G_{\rho}(z)}
 + \frac{G''_{\rho}(z)}{2\,G'_{\rho}(z)}
= \frac{G'_{\rho}(z)}{G_{\rho}(z)}
  \left(\frac{1-g(z)}{2\,g^2(z)}-1\right)
 + \frac{G''_{\rho}(z)}{2\,G'_{\rho}(z)} .
\end{equation} In \cite[Section~4]{RS07} it is shown that the Cauchy transform of $\rho=\mu\boxtimes\nu$ satisfies the algebraic equation
\[
F\big(z,G_\rho(z)\big)=0,
\qquad F(z,y):=z^2y^4-zy+1 .
\]
Differentiating this identity with respect to $z$ gives
\[
\frac{\partial F}{\partial z}\big(z,G_\rho(z)\big)
+\frac{\partial F}{\partial y}\big(z,G_\rho(z)\big)\,G'_\rho(z)=0 .
\]
Since $F_z=2zy^4-y$ and $F_y=4z^2y^3-z$, this yields
\begin{equation}\label{deri}
\frac{G'_{\rho}(z)}{G_{\rho}(z)}
=\frac{1-2zG_{\rho}^3(z)}{z\,\big(4zG_{\rho}^3(z)-1\big)} .
\end{equation}
After some manipulations, using $F(z,G_\rho(z))=0$ and $1-g^2=G_\rho^2$ we may rewrite the above equation as
\begin{equation}\label{deri1}
\frac{G'_{\rho}}{G_{\rho}}=\frac{\,g\,(1+g) G_{\rho}}{2\,(1-2g)} .
\end{equation} Now, writing $G'_\rho=G_\rho\cdot\frac{G'_\rho}{G_\rho}$ and taking logarithmic derivatives in \eqref{deri}, namely
\[
\frac{G''_{\rho}}{G'_{\rho}}=\frac{d}{dz}\log G'_{\rho}
=\frac{G'_\rho}{G_\rho}+\frac{d}{dz}\log\frac{1-2zG_{\rho}^3}{z\,(4zG_{\rho}^3-1)},
\]
and using \eqref{deri} once more to eliminate $G'_\rho$, a straightforward computation gives
\[
\frac{G''_{\rho}}{G'_{\rho}}
=\frac{2\big(24z^3G_{\rho}^9-24z^2G_{\rho}^6+7zG_{\rho}^3-1\big)}
      {z\,\big(1-2zG_{\rho}^3\big)\big(4zG_{\rho}^3-1\big)^2}
=G_{\rho}\,
\frac{(1+g)\big(-1-g+6g^2-6g^3\big)}
     {2\,g\,(1-2g)^2}.
\] 
Substituting this and \eqref{deri1} in \eqref{examff}, we get
\begin{equation*}
\begin{aligned}
G_{\rho'}
&= \frac{G_\rho\,(1+g)}{4g\,(1-2g)^2}
   \Big[(1+g)(1-2g)^2 + \big(-1-g+6g^2-6g^3\big)\Big] \\
&= -\,\frac{G_\rho\,(1-g^2)(2-g)}{2\,(1-2g)^2}.
\end{aligned}
\end{equation*}
Finally, since $1-g^2=4G_\rho^2$, we conclude that \eqref{eq:sec-ex} holds.
\end{exam}

In the next set of examples, we assume that $\mfq=(q_d)_{d\geq t}$ is a finite rank perturbation of the polynomial $(x-\alpha)^d$. To be precise, $q_d$ is as follows, \[
q_d(x)=(x-\alpha)^{d-t}(x-\beta_1)(x-\beta_2)\cdots(x-\beta_t).
\] The infinitesimal distribution of this sequence is given by, \cite[Example 5.3]{ADS25},
\[
(\nu,\nu')=\Big(\delta_{\alpha},\; -t\,\delta_{\alpha}+\sum_{k=1}^{t}\delta_{\beta_k}\Big).
\]
Firstly, since 
$R_{\nu}(z)=\alpha,$ by the statement $(e)$ in Proposition \ref{p2}, it follows that
$$C_1(z)=\alpha \quad \text{and} \quad C_2(z)=R_{\mu}(\alpha z).$$ Therefore, from Theorem \ref{ffp3}, we get
\begin{align} \label{finiterank}\widehat{R}_{\mfp\boxtimes\mfq}(z)&=\alpha \,\widehat{R}_{\mfp}(\alpha z)+ (z R_{\mu}(\alpha z))^{'}\widehat{R}_{\mfq}(z R_{\mu}(\alpha z))
\end{align} 
We will consider the two main cases, $\alpha=1$ and $\alpha=0$. The case for general $\alpha\neq0$ may be derived by a rescaling of the case $\alpha=1$.

\begin{exam}  [Multiplication with finite rank polynomial ($\alpha=0$)]
Assume that $\mfq=(q_d)_{d\geq t}$ is as follows \[
q_d(x)=(x)^{d-t}(x-\beta_1)(x-\beta_2)\cdots(x-\beta_t).
\] This corresponds to $\alpha=0$ in the above considerations. In this case, we see that $\alpha \,\widehat{R}_{\mfp}(\alpha z)=0$ and $R_\mu(\alpha z)=R_\mu(0)=\kappa_{1}^{\mu}$, so that
$$\widehat{R}_{\mfp\boxtimes\mfq}(z)= \kappa_{1}^{\mu}\widehat{R}_{\mfq}(\kappa_{1}^{\mu} z).$$ Moreover, $\rho=\delta_0\boxtimes \nu=0$, and thus $G_\rho(z)=\frac{1}{z}$ and $H_\rho=0$. Therefore, \eqref{qis0} yields
\[
G_{\rho'}(z)=\kappa_{1}^{\mu}\cdot\widehat{R}_{\mfq}\left(\frac{\kappa_{1}^{\mu}}{z}\right)\cdot\frac{1}{z^2}
=
\sum_{k=1}^{t}
\frac{ \kappa_{1}^{\mu}}{z-\beta_k \kappa_{1}^{\mu}}- \frac{t\kappa_{1}^{\mu}}{z}.
\]
or equivalently,
$\rho'=-t\,\kappa_{1}^{\mu}\,\delta_0+\kappa_{1}^{\mu}\sum_{k=1}^{t}\delta_{\kappa_{1}^{\mu}\beta_k}$. This result is aligned with Example 5.16 from \cite{ADS25}.

\end{exam}

\begin{exam}[Multiplication with finite rank polynomial ($\alpha=0$)]
Assume now that $\mfq=(q_d)_{d\geq t}$ is given by
$$
q_d(x)=(x-1)^{d-t}(x-\beta_1)(x-\beta_2)\cdots(x-\beta_t),
$$
which corresponds to the case $\alpha=1$ in the discussion above. In particular,
$\rho=\mu\boxtimes\delta_1=\mu$. By \eqref{finiterank}, we have
$$
\widehat{R}_{\mfp\boxtimes\mfq}(z)
=\widehat{R}_{\mfp}(z)+\big(zR_{\mu}(z)\big)'\,\widehat{R}_{\mfq}\big(zR_{\mu}(z)\big).
$$
Substituting this into \eqref{Grho}, we obtain
\begin{align*}
G_{\rho'}(z)
&= -\widehat{R}_{\mfp\boxtimes\mfq}\big(G_{\mu}(z)\big)\,G'_{\mu}(z)-\gdos_{\mu}(z)\\
&= -\widehat{R}_{\mfp}\big(G_{\mu}(z)\big)\,G'_{\mu}(z)-\gdos_{\mu}(z)\\
&\qquad -\Big[G_\mu(z)\,G_\mu'(z)\,R'_\mu\big(G_\mu(z)\big)+G_\mu'(z)\,R_\mu\big(G_\mu(z)\big)\Big]
\,\widehat{R}_{\mfq}\Big(G_\mu(z)\,R_\mu\big(G_\mu(z)\big)\Big)\\
&= G_{\mu'}(z)-\Big(G_\mu(z)\,R_\mu\big(G_\mu(z)\big)\Big)'
\,\widehat{R}_{\mfq}\Big(G_\mu(z)\,R_\mu\big(G_\mu(z)\big)\Big)\\
&= G_{\mu'}(z)-\big(zG_\mu(z)\big)'\,\widehat{R}_{\mfq}\big(zG_\mu(z)-1\big),
\end{align*}
where we used in the last equality we used the relation $R_\mu\big(G_\mu(z)\big)=z-\frac{1}{G_\mu(z)}$, which gives
$G_\mu(z)\,R_\mu\big(G_\mu(z)\big)=zG_\mu(z)-1.$

Moreover, by \cite[Example~5.3]{ADS25},
\[
  \hat{R}_{\mathfrak{q}}(z) = \sum_{k=1}^{t} \frac{\beta_k - 1}{1 - (\beta_k - 1) z}.
\]
Substituting this expression into the previous identity yields
\begin{equation}\label{Gro1}
    G_{\rho'}(z) = G_{\mu'}(z) + \big(zG_\mu(z)\big)^{'}\sum_{k=1}^{t} \frac{1-\beta_k }{\beta_k + (1-\beta_k) zG_{\mu}(z)}.  
\end{equation} In particular, if $\beta_k = 0$ for every $k \in [t]$, this formula reduces to
\[
  G_{\rho'}(z) 
  = G_{\mu'}(z) + t \left( \frac{G_{\mu}'(z)}{G_{\mu}(z)} + \frac{1}{z} \right),
\]
in agreement with \cite[Example~5.18]{ADS25}.
\end{exam}

\begin{exam}  [Multiplication of finite rank polynomials ($\alpha=1$)]
Assume that $\mfp=(p_d)_{d\geq s}$ and $\mfq=(q_d)_{d\geq t}$ are as follows
\[
p_d(x)=(x-1)^{d-s}(x-\alpha_1)(x-\alpha_2)\cdots(x-\alpha_s),
\qquad 
q_d(x)=(x-1)^{d-t}(x-\beta_1)(x-\beta_2)\cdots(x-\beta_t).
\]
This is a particular case of the previous example. However, in this case we have more information, namely, $G_{\mu}(z)=\frac{1}{z-1}$. Thus, according to \eqref{Gro1}
\begin{align*}
G_{\rho'}(z)
&= G_{\mu'}(z) + \Big(\frac{z}{z-1}\Big)' \sum_{k=1}^{t} \frac{1-\beta_k}{\beta_k + (1-\beta_k)\frac{z}{z-1}} \\
&= G_{\mu'}(z) + \sum_{k=1}^{t} \frac{\beta_k-1}{(z-\beta_k)(z-1)} \\
&= G_{\mu'}(z) + G_{\nu'}(z) \\
&= -\frac{s+t}{z-1} + \sum_{k=1}^{s}\frac{1}{z-\alpha_k} + \sum_{k=1}^{t}\frac{1}{z-\beta_k},
\end{align*} or equivalently,
$\rho'=-(s+t)\,\delta_1+\sum_{k=1}^{s}\delta_{\alpha_k}+\sum_{k=1}^{t}\delta_{\beta_k}$. This result is aligned with Example 5.17 from \cite{ADS25}.

\end{exam}

\bibliography{Bib}
\bibliographystyle{acm}

\end{document}